\documentclass{article}
\usepackage{graphicx} 
\usepackage{tikz,amsmath,pgfplots,siunitx,bm}
\usepackage{authblk}
\usepackage{multirow}
\usepackage{comment}
\usepackage{subfigure}
\usepackage[makeroom]{cancel}
\usepackage[framed,numbered,useliterate]{mcode}
\newcommand{\bx}{\bm{x}}
 
\usepackage{soul}

\usepackage{listings}
\usepackage{pythonhighlight}
\lstnewenvironment{pythonLines}[1][]{\lstset{style=mypython,numbers=left}}{}

\newtheorem{remark}{Remark}

\title{A simple GPU implementation of spectral-element methods for solving 3D Poisson type equations on rectangular domains and its applications} 

\author[1]{Xinyu Liu\thanks{liu1957@purdue.edu}}
\author[2]{Jie Shen\thanks{jshen@eitech.edu.cn}}
\author[1]{Xiangxiong Zhang\thanks{zhan1966@purdue.edu}}

\affil[1]{Department of Mathematics, Purdue University, West Lafayette, IN 47906, US} 
\affil[2]{Eastern Institute of Technology, Ningbo, Zhejiang 315200, P. R. China} 
\date{}
\pgfplotsset{compat=1.18} 
\begin{document}

\maketitle
\begin{abstract}
It is well known since 1960s that by exploring the tensor product structure  of the discrete Laplacian  on Cartesian meshes, one can develop a simple direct Poisson solver with an $\mathcal O(N^{\frac{d+1}d})$ complexity in $d$-dimension, where $N$ is the number of the total unknowns.   
The GPU acceleration of numerically solving PDEs has been explored  successfully  around fifteen years ago and become more and more popular in the past decade, driven  by significant advancement in both hardware and software technologies, especially in the recent few years.   We present in this paper a simple but extremely fast MATLAB implementation on a modern GPU,  which can be easily reproduced, for solving 3D Poisson type equations using a spectral-element method. In particular, it costs less than one second on a Nvidia A100 for solving a Poisson equation with one billion degree of freedoms. We also present applications of this fast solver to solve a linear (time-independent) Schr\"odinger equation and a nonlinear (time-dependent) Cahn-Hilliard equation.
\end{abstract}

\section{Introduction}

It is well known that the tensor product structure of  the discrete Laplacian on Cartesian meshes can be used to invert the Laplacian since 1960s \cite{lynch1964tensor}. This approach   has been particularly popular for spectral and spectral-element methods  \cite{HAIDVOGEL1979167, PATERA1986474, shen1994, kwan2007efficient,ChenShen2012}. 
In fact, this method can be used for any discrete Laplacian  on a Cartesian mesh.
In this paper, as an example, we focus on the $Q^k$ spectral-element method, which is equivalent to the classical continuous finite element method with Lagrangian $Q^k$ basis implemented with the  $(k+1)$-point Gauss--Lobatto quadrature \cite{maday1990optimal, li2022accuracy}. 
Tensor based  solvers naturally fit the design of graphic processing units (GPUs). The earliest successful attempts to accelerate the computation of high order accurate methods in scientific computing communities include nodal discontinuous Galerkin method \cite{klockner2009nodal} almost fifteen years ago. These pioneering efforts of GPU acceleration of high order methods, or even those ones published later such as \cite{chen2013gpu} in 2013, often rely on intensive coding.

In recent years, the surge in computational demands from machine learning and neural network based approaches has led to the evolution of modern GPUs. Correspondingly, software technologies have advanced considerably, streamlining the utilization of
GPU computing.
The landscape of both hardware and software has dramatically transformed, differing substantially from what existed a decade or even just two years ago.

In this paper, we present a straightforward yet robust implementation of accelerating the spectral-element method for three-dimensional discrete Laplacian 
 on modern GPUs. In particular, for a total number of degree of freedoms as large as $1000^3$, the inversion of the 3D Laplacian using an arbitrarily high order $Q^k$ spectral-element method, takes no more than one second on one Nvidia A100 GPU card with 80G memory. 
While this impressive computational speed is naturally contingent on the hardware, it is noteworthy that our approach is grounded in a minimalist MATLAB implementation, ensuring ease of replication. In the Appendix, we give a full MATLAB code for solving a 3D Poisson equation on a rectangular domain using $Q^k$ spectral-element method. 

We remark that a similar simple implementation  on GPUs can also be achieved using Python using the Python package JAX, which however provides better performance than MATLAB only for single precision computation \textcolor{blue}{under TensorFloat-32 precision format}.  Our numerical tests and comparison on one Nvidia A100 GPU card suggest that MATLAB implementation performs better than Python for double precision computation for large problems like one billion DoFs.

We emphasize that the ability of solving Poisson type equation fast can play an important role in many fields of science and engineering. In fact, a large class of time dependent  nonlinear systems, after a suitable implicit-explicit (IMEX) time discretization,   often reduces to solving  Poisson type equations at each time step (see, for instance, \cite{shen2019new}). Therefore, having a simple, accurate  and very fast solver for Poisson type equations can lead to very efficient numerical algorithms  on modern GPUs for these nonlinear systems which include, e.g., Allen-Cahn and Chan-Hillard equations and related phase-field models \cite{shen2019new}, nonlinear Schr\"odinger equations, Navier-Stokes equations and related hydro-dynamic equations through a decoupled (projection, pressure correction etc.) approach \cite{guermond2006overview}. In particular, by using the code provided in the Appendix, one can build, with a relatively easy effort,  very efficient numerical   solvers  on modern GPUs for these time dependent complex nonlinear systems.

The rest of the paper is organized as follows. In Section \ref{sec:tensor}, we give the  implementation details for 3D problems,  which is robust for very high order elements with the computation approach in Section \ref{remark:robustness}.
In Section \ref{sec:test}, we demonstrate  the good performance of this simple implementation for equations including the Poisson equation,  a variable coefficient elliptic problem solved by the preconditioned conjugate gradient descent using Laplacian as a preconditioner, 
as well as a Cahn--Hilliard equation.  Although our focus in this paper remains on the spectral-element method for these particular equations, similar results can be obtained  for other problems with the same tensor product structure,  e.g., finite difference schemes in implementing the matrix exponential in the exponential time differencing \cite{du2019maximum} and spectral fractional Laplacian \cite{chen2020efficient}. 
We also compare it with a similar simple implementation in Python, and the numerical results suggest that MATLAB implementation is better than Python for double precision computation on A100. 
Some concluding remarks are given in Section \ref{sec:remark}.

\section{A spectral-element method for Poisson type equations}
\label{sec:tensor}
To fix the idea, we describe the implementation details for solving the Poisson type equation
\begin{equation}\label{pde}
    \alpha u - \Delta u = f,
\end{equation}
with a constant coefficient $\alpha>0$ and homogeneous Neumann boundary conditions on a  rectangular domain $\Omega$.  In this paper,  we only consider  the spectral-element method with continuous piecewise $Q^k$ polynomial basis on uniform rectangular meshes, and all integrals are approximated by $(k+1)$-point Gauss--Lobatto quadrature \cite{li2020superconvergence}.

\subsection{The spectral-element method in two dimensions}
\begin{figure}[htbp]
 \subfigure[A $3\times 2$ mesh for $Q^2$ element and the $3\times 3$ Gauss-Lobatto quadrature. ]{\includegraphics[scale=0.8]{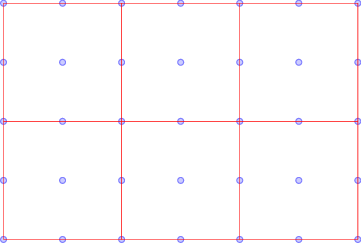}}
 \hspace{.6in}
 \subfigure[All the  quadrature points for $Q^2$ element. ]{\includegraphics[scale=0.8]{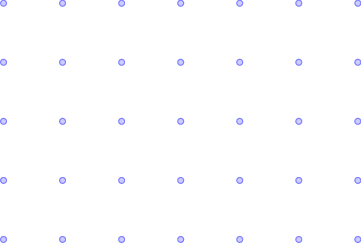}}\\
       \subfigure[A $2\times 2$ mesh for $Q^3$ element and the $4\times 4$ Gauss-Lobatto quadrature. ]{\includegraphics[scale=0.8]{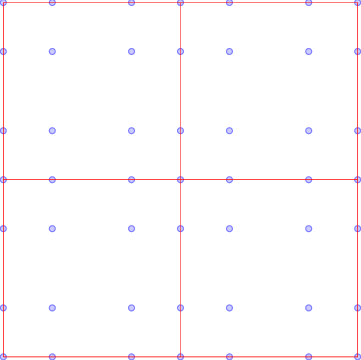} }
 \hspace{.6in}
 \subfigure[All the  quadrature points for $Q^3$ element. ]{\includegraphics[scale=0.8]{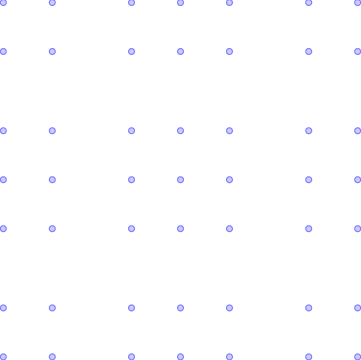}}
\caption{An illustration of Lagrangian $Q^k$ element and the $(k+1)\times(k+1)$ Gauss-Lobatto quadrature. }
\label{mesh-Q2}
 \end{figure}

We first consider the two dimensional case. 
As shown in Figure \ref{mesh-Q2}, 
such quadrature points naturally define all degree of freedoms since a single variable polynomial of degree $k$ is uniquely determined by its values at $k+1$ points.

On a rectangular mesh for $Q^k$ basis as shown in  Figure \ref{mesh-Q2} (a) or (c), let $(x_i, y_j)$ ($i=1,\cdots, N_x; j=1,\cdots N_y$) denote all the points in Figure \ref{mesh-Q2} (b) or (d). We consider the following  finite element space of continuous piecewise $Q^k$ polynomials:
\begin{equation*}
    V^{h} = Span\{\phi_{i}(x)\phi_{j}(y), 1\leq i\leq N_x, 1\leq j\leq N_y\},
\end{equation*}
where $\phi_i(x)$ ($i=1,\cdots, N_x$) denotes the $i$-th Lagrangian interpolation polynomial of degree $k$ in the $x$ direction as shown in  Figure \ref{mesh-Q2} (b) or (d). 

Then, the $Q^k$ spectral-element method for solving \eqref{pde} on a rectangular domain $\Omega$ is to seek $u_h\in V^h$ satisfying
\begin{equation}
    \label{fem}
\alpha\langle u_h,  v_h \rangle+\langle \nabla u_h, \nabla v_h \rangle =
\langle f, v_h\rangle,\quad \forall v_h\in V^h,
\end{equation}
where $\langle f, g\rangle$ denotes the approximation to the integral $\iint_\Omega f(x,y) g(x,y) \, dx dy$ by the $(k+1)\times (k+1)$ Gauss-Lobatto quadrature rule in each cell as shown in Figure \ref{mesh-Q2}.

The numerical solution $u_h(x,y)\in V^h$  can be expressed by the basis as 
\begin{equation*}
    u_h(x,y) = \sum_{i=1}^{N_x}\sum_{j=1}^{N_y} u_{i,j}\phi_{i}(x)\phi_{j}(y),
\end{equation*}
where the coefficients $u_{i,j}=u_h(x_i,y_j)$ since $\phi_i(x)$ and $\phi_j(y)$ are chosen as the Lagrangian interpolation polynomials at $x_i$ and $y_j$.

Next, we define the one-dimensional stiffness matrix and the mass matrix as follows. The stiffness matrix $S_x$ is a matrix of size $N_x\times N_x$ with $(i,j)$-th entry being $\langle \phi'_i(x), \phi'_j(x)\rangle$.
The mass matrix $M_x$ is a matrix of size $N_x\times N_x$ with $(i,j)$-th entry being $\langle  \phi_i(x),  \phi_j(x)\rangle$. 
The matrices $S_y$ and $M_y$ are similarly defined, with a size $N_y\times N_y$. Since the basis polynomials $\phi_i(\cdot)$ are Lagrangian interpolants at the quadrature points, the mass matrices $M_x$ and $M_y$ are diagonal. 

\begin{remark}
 Instead of the Lagrangian basis, we can also use the modal basis $L_j(x)-L_{j+2}(x)$ in each cell, where $L_k(x)$ is the Legendre polynomial of degree $k$,  as the interior basis functions  and the piecewise linear hat function at the intersecting points of two subintervals. This leads to sparse mass and stiffness matrices \cite{kwan2007efficient,ChenShen2012}.
\end{remark}
 
Let $U$ be a matrix of size $N_x\times N_y$ with $u_{i,j}$ being its $(i,j)$-entry. Then  $u_h(x,y)$ can be equivalently represented by the matrix $U$.
Similarly, let $F$ be a matrix of size $N_x\times N_y$ with $f(x_i, y_j)$ being its $(i,j)$-entry. 
Then the scheme  \eqref{fem} can be equivalently written as
\begin{equation}\label{weak2D}
    \alpha M_{x}UM^T_{y} + S_{x}UM^T_{y} + M_{x}US^T_{y} = M_{x} F M^T_{y}.
\end{equation}
Detailed derivations of \eqref{weak2D} can be found in \cite{li2020superconvergence,shen2022discrete, hu2023positivity}.

\subsection{Inversion by eigenvalue decomposition}

For any matrix $X$ of size $m\times n$, define a vectorization operation $vec(\cdot)$, and let $vec(X)$ be the vector of size $mn$ obtained by reshaping all entries of $X$ into a column vector in a column by column order. Then for any two matrices $A_1, A_2$ of proper sizes, it satisfies
\begin{equation}\label{basic2Dvec}
    vec(A_1XA_2^T) = \left(A_2\otimes A_1\right)vec(X),
\end{equation}
where $\otimes$ denotes the Kronecker product. With \eqref{basic2Dvec}, the $Q^k$ spectral-element method \eqref{weak2D} is also equivalent to 
\begin{equation}\label{matvec2D}
    (\alpha M_{y}\otimes M_{x} + {M_{y}\otimes S_{x} + S_{y}\otimes M_{x}})vec(U) = (M_{y}\otimes M_{x})vec(F).
\end{equation}

{  We remark that many other types of boundary conditions on a cubic domain can be written in the form of \eqref{weak2D} or \eqref{matvec2D} include periodic boundaries and the simple homogeneous Dirichlet boundary treatment described   in \cite{li2020superconvergence}. For  non-homogeneous Neumann and Dirichlet boundary conditions, the boundary conditions will be added to the right hand side of \eqref{weak2D} or \eqref{matvec2D}, and the left hand side of \eqref{weak2D} or \eqref{matvec2D} stays the same, i.e., the matrix to be inverted stays the same. The left hand side of   \eqref{matvec2D} would no longer hold if the boundary value problem does not have a tensor product structure, e.g., on a curved domain.}

The linear system \eqref{weak2D}, or equivalently  \eqref{matvec2D}, can be solved by the following well-known method by using eigenvalue decomposition for only small matrices such as $S$ and $M$. For convenience, we only consider the simplified equivalent system 
\begin{equation}\label{weak2D-2}
    \alpha U  + H_{x}U  +  UH^T_{y} =  F,
\end{equation}
or 
\begin{equation}\label{matvec2D-3}
    (\alpha  + {I_{y}\otimes H_{x} + H_{y}\otimes I_{x}})vec(U) = vec(F),
\end{equation}
where $I$ is the identity matrix and $H=M^{-1}S$.

First, solve a generalized eigenvalue problem for small matrices $S$ and $M$, i.e., finding eigenvalues $\lambda_i$ and eigenvectors $\mathbf v_i$ satisfying
\begin{equation}
    S\mathbf v_i =\lambda_i M \mathbf v_i. \label{generalized-eig}
\end{equation}

Regardless of what kind of basis functions is used in a spectral-element method, the variational form of \eqref{fem} ensures the symmetry of $S$ and $M$, thus a complete set of eigenvectors exists for \eqref{generalized-eig}.
Let $\Lambda$ be a diagonal matrix with all eigenvalues $\lambda_i$ being diagonal entries, and let $T$ be the matrix with all corresponding eigenvectors $\mathbf v_i$ as its columns. Then 
\[ST=M T \Lambda\Rightarrow H=M^{-1}S=T\Lambda T^{-1}. \]
Thus \eqref{matvec2D-3} becomes
\begin{equation*} 
    [\alpha  + {(T_y I_y T_y^{-1})\otimes (T_x \Lambda_x T_x^{-1}) + (T_y\Lambda_y T_y^{-1})\otimes (T_x I_x T_x^{-1})}]vec(U) = vec(F),
\end{equation*}
which is equivalent to
\begin{equation}\label{matvec2D-4}
    (T_y\otimes T_x)[\alpha  + {I_y \otimes   \Lambda_x + \Lambda_y  \otimes  I_x}] (T_y^{-1}\otimes T_x^{-1})vec(U) = vec(F).
\end{equation}
Notice that $\alpha  +  I_y \otimes   \Lambda_x +\Lambda_y  \otimes  I_x$ is a diagonal matrix,   thus its inverse is simple to compute. 
Let $\Lambda2D$ be a matrix of size $N_x\times N_y$ with its $(i,j)$ entry being equal to $\alpha+(\Lambda_x)_{i,i}+(\Lambda_y)_{j,j}$, then  \eqref{matvec2D-4} can be solved by  
\begin{equation}
    U= T_x[(T_x^{-1}F T_y^{-T})./\Lambda2D]T_y^{T}, 
    \label{tensormethod-2D}
\end{equation}
where $./$ denotes the entrywise division between two matrices.

\subsection{Robust computation of the generalized eigenvalue problem}
\label{remark:robustness}
In our spectral-element implementation, $M$ is diagonal. So we can consider an eigenvalue problem instead of the generalized eigenvalue problem \eqref{generalized-eig}.
A numerically robust method, especially for very high order polynomial basis, is   to solve the following symmetric eigenvalue problem. Let
\begin{equation*}
    H = M^{-1}S = M^{-1/2} (M^{-1/2}SM^{-1/2}) M^{1/2}.
\end{equation*}
Since $S_1 = M^{-1/2}SM^{-1/2}$ is real and symmetric, we can first find its eigenvalue decomposition as $S_1 = Q\Lambda Q^T$ where $\Lambda$ is a diagonal matrix and $Q$ is an orthogonal matrix. Then, we have 
\begin{equation*}
    H = M^{-1/2} (Q\Lambda Q^T) M^{1/2} = T\Lambda T^{-1},
\end{equation*}
with $T = M^{-1/2}Q$ and $T^{-1} = Q^T M^{1/2}$. In Section \ref{sec:robustness}, we will show numerical tests validating the robustness of this implementation for very high order elements.

\subsection{Implementation for the three-dimensional case}
\label{sec-3d}
On a three dimensional rectangualar mesh, any continuous piecewise $Q^k$ polynomial $u_h$
can be uniquely represented by a 3D array $U$ of size $N_x\times N_y\times N_z$ with $(i,j,k)$-th entry denoting the point value $u_h(x_i,y_j,z_k)$, where $(x_i, y_j, z_k)$ ($i=1,\cdots, N_x; j=1,\cdots, N_y; k=1,\cdots, N_z$) denotes all the quadrature points. 

For a 3D array $U$, we define a page as the matrix obtained by fixing the last index of $U$. Namely, $U(\,:\,,\,:\,,k)$ for any fixed $k$ is a page of $U$. For a matrix $U(\,:\,,\,:\,,k)$ of size $N_x\times N_y$, recall that $vec(U(\,:\,,\,:\,,k))$ is a column vector of size $N_xN_y$.
We define $\hat U$ as the following matrix of size $N_xN_y\times N_z$ obtained by reshaping $U$:
\[\hat U=\begin{bmatrix}
    vec(U(\,:\,,\,:\,,1)) & vec(U(\,:\,,\,:\,,2)) & \cdots & vec(U(\,:\,,\,:\,,N_z)) 
\end{bmatrix}.\]
Then we define $vec(U)$ as the vector of size $N_xN_yN_z\times 1$ by reshaping $\hat U$ in a column by column order.

With the notation above, it is straightforward to verify that
\begin{equation}
    \label{vec-array}
     (A_3^T\otimes A_2^T\otimes A_1)vec(U) =vec((A_2^T\otimes A_1)\hat U A_3).
\end{equation}

Next, we consider how to implement the matrix vector multiplication in \eqref{vec-array}  without reshaping the 3D arrays.
Let $Y$ be a 3D array of size $N_x\times N_y\times N_z$ defined by
\begin{equation}
    \label{tensorprod}
    vec(Y)=(A_3^T\otimes A_2^T\otimes A_1)vec(U).
\end{equation}
With the simple property \eqref{vec-array}, in our numerical tests, we find that the following simple implementation of \eqref{tensorprod} in MATLAB 2023 is efficient using two functions {\it tensorprod} and {\it pagemtimes}:
\begin{lstlisting}
% Computing a 3D array Y of the same size as U defined above
Y = tensorprod(U,A3,3,1);
Y = pagemtimes(Y,A2);
Y = squeeze(tensorprod(A1,Y,2,1));
\end{lstlisting}

For the three-dimensional case, for simplicity, we consider the equation \eqref{pde}  with $\alpha=0$. With similar notation as in the two-dimensional case, the matrix form of the $Q^k$ spectral-element method \eqref{fem} can be given as
\begin{eqnarray*}
    \left(M_{z}\otimes M_{y}\otimes S_{x} + M_{z}\otimes S_{y}\otimes M_{x} + S_{z}\otimes M_{y}\otimes M_{x} \right)vec(U) \\
    =\left(M_{z}\otimes M_{y}\otimes M_{x} \right)vec(F),
\end{eqnarray*}
or equivalently,
\begin{equation}\label{matvec3D}
    \left(I_{z}\otimes I_{y}\otimes H_{x} + I_{z}\otimes H_{y}\otimes I_{x} + H_{z}\otimes I_{y}\otimes I_{x} \right)vec(U) =vec(F),
\end{equation}
where $U$ is a 3D array with $(i,j,k)$-th entry denoting the point value $u_h(x_i,y_j,z_k)$, and $F$ is a 3D array with $(i,j,k)$-th entry denoting the point value $f(x_i,y_j,z_k)$.

With the eigenvalue decomposition $H=M^{-1}S=T\Lambda T^{-1}$,
similar to the derivation of \eqref{matvec2D-4},
 the equation
\eqref{matvec3D} is equivalent to 
\begin{equation}
    \label{3D-Poisson}
\resizebox{.99\hsize}{!}{$ (T_z \otimes T_y \otimes T_x) \left( I_{z}\otimes I_{y}\otimes \Lambda_{x} + I_{z}\otimes \Lambda_{y}\otimes I_{x} + \Lambda_{z}\otimes I_{y}\otimes I_{x} \right)(T_z^{-1}\otimes T_y^{-1}\otimes T_x^{-1}) vec(U) =vec(F)$}.
\end{equation}
Define a 3D array $\Lambda3D$ with its $(i,j,k)$-th entry being equal to $(\Lambda_x)_{i,i}+(\Lambda_y)_{j,j}+(\Lambda_z)_{k,k}$, then 
\eqref{3D-Poisson}
can be implemented efficiently as the following in MATLAB:
\begin{table}[htbp]
    \centering
  \begin{lstlisting}
% Simple and efficient implementation of (14)
% TInv denotes the inverse matrix of T
U = tensorprod(F,TzInv',3,1);
U = pagemtimes(U,TyInv');
U = squeeze(tensorprod(TxInv,U,2,1));
U = U./Lambda3D;
U = tensorprod(U,Tz',3,1);
U = pagemtimes(U,Ty');
U = squeeze(tensorprod(Tx,U,2,1));
\end{lstlisting}
    \caption{\it The MATLAB 2023 script of implementing \eqref{3D-Poisson} on both CPU and GPU.}
    \label{tab:3dcode}
\end{table}

\section{Numerical tests}
\label{sec:test}

In this section, we report the performance of the simple MATLAB implementation in Table \ref{tab:3dcode}. In particular, a demonstration code is provided in the Appendix. The performance and speed-up are of course dependent on the hardwares. We test our code on the following three devices:
\begin{enumerate}
    \item CPU: Intel  i7-12700  2.10 GHz (12-core) with 16G memory;
    \item GPU: Quadro RTX 8000 (48G memory); 
    \item GPU: Nvidia A100 (80G memory).
\end{enumerate}

    In MATLAB 2023, for computation on either CPU or GPU, the code for implementing \eqref{3D-Poisson} is the same as in Table \ref{tab:3dcode}. On the other hand, matrices like $T_x, T_y, T_z$ and arrays like $F$ and $\Lambda3D$ must be loaded to GPU memory before performing the GPU computation, see the full code in the Appendix. We define the process of loading matrices and arrays $T_x, T_y, T_z, F, \Lambda3D$ as the offline step since it is preparatory, and undertaken only once, regardless of how many times the Laplacian needs to be inverted. We define the step in Table \ref{tab:3dcode} as the online computation step. {\bf All the computational time reported in this section are online computational time, i.e., we do not count the offline preparational time.} 
     
\subsection{Accuracy tests}
\label{sec:poisson-accuracy}
We list a few accuracy tests to show that 
the scheme implemented is indeed high order accurate. 
In particular, the $Q^k$ ($k\geq 2$) spectral-element method is $(k+2)$-th order accurate for smooth solutions when measuring the $\ell^2$ error in function values for solving second order PDEs, which has been rigorously proven recently 
in \cite{li2020superconvergence, li2022accuracy}. 

We consider the Poisson type equation \eqref{pde}
with $\alpha = 1$ in domain $\Omega = [-1,1]^3$. For Dirichlet boundary conditions, we test a smooth exact solution
\begin{equation*}
    u_{D}^{*}(\bx) = \sin(\pi x)\sin(2\pi y)\sin(3\pi z) + (x-x^3)(y^2-y^4)(1-z^2).
\end{equation*}
For Neumann boundary conditions, we test a smooth exact solution
\begin{equation*}
    u_{N}^{*}(\bx) = \cos(\pi x)\cos(2\pi y)\cos(3\pi z) + (1-x^2)^3(1-y^2)^2(1-z^2)^4.
\end{equation*}
The results of $Q^5$ and $Q^6$ spectral-element methods are listed in Table \ref{tab1:accuracy test dirichlet/neumann}.  
\begin{table}[htbp]
    \centering
    \begin{tabular}{|c|c c c|c c c|}
    \hline
    \multicolumn{7}{|c|}{$Q^5$ spectral-element method (SEM)}\\
    \hline
    \multirow{2}{*}{FEM Mesh} &  \multicolumn{3}{c|}{Dirichlet boundary} &  \multicolumn{3}{c|}{Neumann boundary}\\
     & {Total DoFs}  & $\ell^2$ error & order & {Total DoFs}  & $\ell^2$ error & order\\
    \hline
    $2^3$ &$9^3$ &  $2.27$E-$1$ & - & $11^3$ & $4.76$E-$1$ & - \\
    \hline
    $4^3$ & $19^3$ &  $3.91$E-$3$ & $5.86$ & $21^3$ & $5.49$E-$3$ & $6.44$\\
    \hline
    $8^3$ & $39^3$ &  $4.12$E-$5$ & $6.57$ & $41^3$ & $4.32$E-$5$ & $6.99$\\
    \hline
    $16^3$ & $79^3$ & $3.34$E-$7$ & $6.95$ & $81^3$ & $3.42$E-$7$ & $6.98$\\
    \hline
    $32^3$ & $159^3$ & $2.63$E-$9$ & $6.99$ & $161^3$ & $2.67$E-$9$ & $7.00$\\
    \hline
    \hline
    \multicolumn{7}{|c|}{$Q^6$ spectral-element method (SEM)}\\
    \hline
    \multirow{2}{*}{FEM Mesh} &  \multicolumn{3}{c|}{Dirichlet boundary} & \multicolumn{3}{c|}{Neumann boundary}\\
     & {Total DoFs}  & $\ell^2$ error & order & {Total DoFs}  & $\ell^2$ error & order\\
    \hline
    $2^3$ & $11^3$ &  $9.68$E-$2$ & - & $13^3$ & $1.18$E-$1$ & - \\
    \hline
    $4^3$ & $23^3$ &  $6.05$E-$4$ & $7.32$ & $25^3$ & $8.42$E-$4$ & $7.13$\\
    \hline
     $8^3$ & $47^3$ & $3.11$E-$6$ & $7.60$ & $49^3$ & $3.24$E-$6$ & $8.02$\\
    \hline
    $16^3$ & $95^3$ & $1.26$E-$8$ & $7.95$ & $97^3$ & $1.28$E-$8$ & $7.98$\\
    \hline
    $32^3$ & $191^3$ & $4.96$E-$11$ & $7.98$ & $193^3$ & $5.09$E-$11$ & $7.98$\\
    \hline
    \end{tabular}
    \caption{\it Accuracy tests for discrete Laplacian for a 3D problem with Dirichlet boundary conditions and a 3D problem with Neumann boundary conditions.}
    \label{tab1:accuracy test dirichlet/neumann}
\end{table}

\subsection{GPU acceleration for solving a Poisson type equation}
\label{sec:poisson-speedup}

\begin{figure}[htbp]
 \subfigure[Comparison on three devices.]{\includegraphics[scale=0.18]{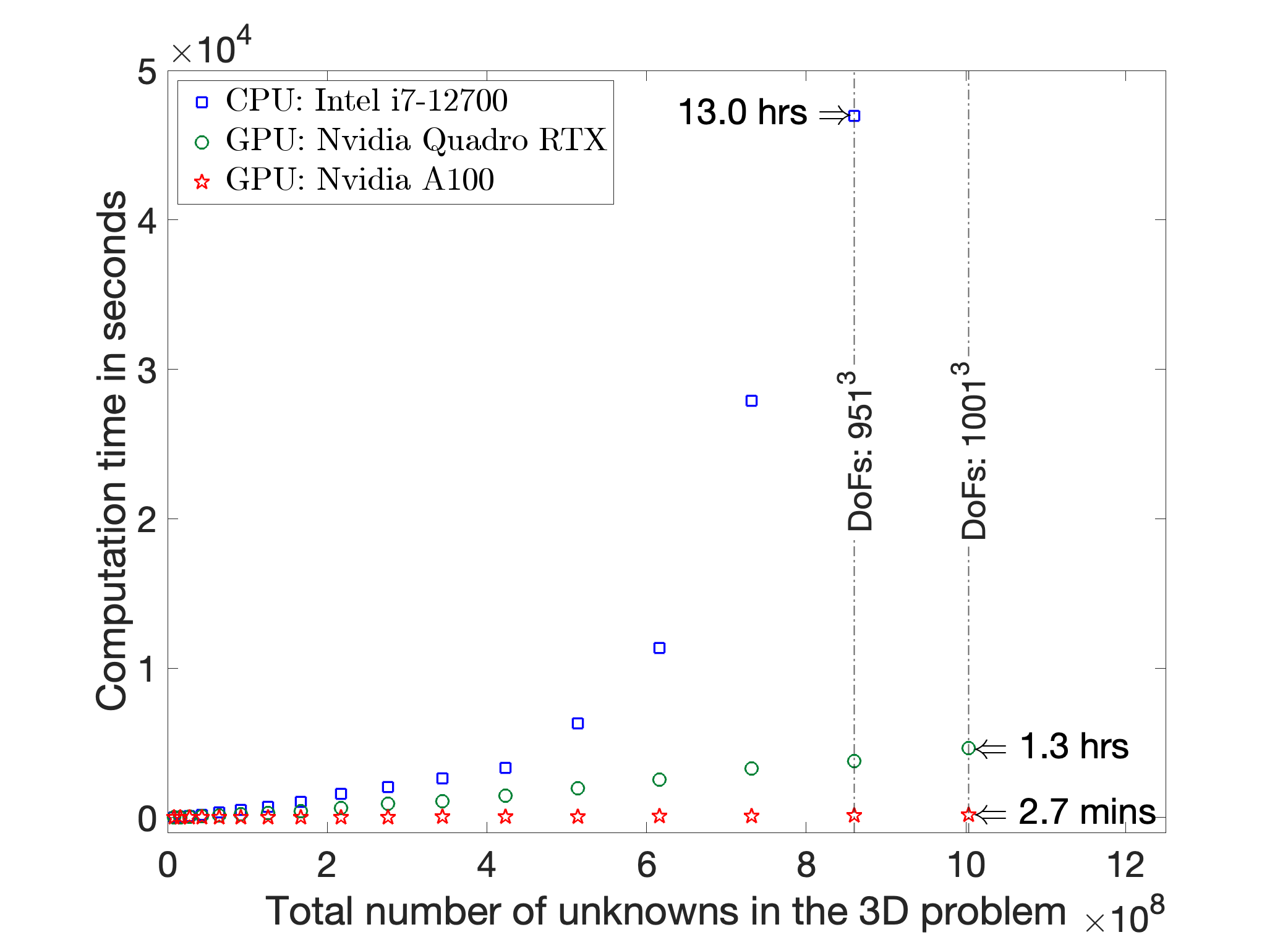}}
 \subfigure[Semilogx plot shows the complexity on all devices are $\mathcal O(N^{\frac43})$ for problems with proper sizes.]{\includegraphics[scale=0.18]{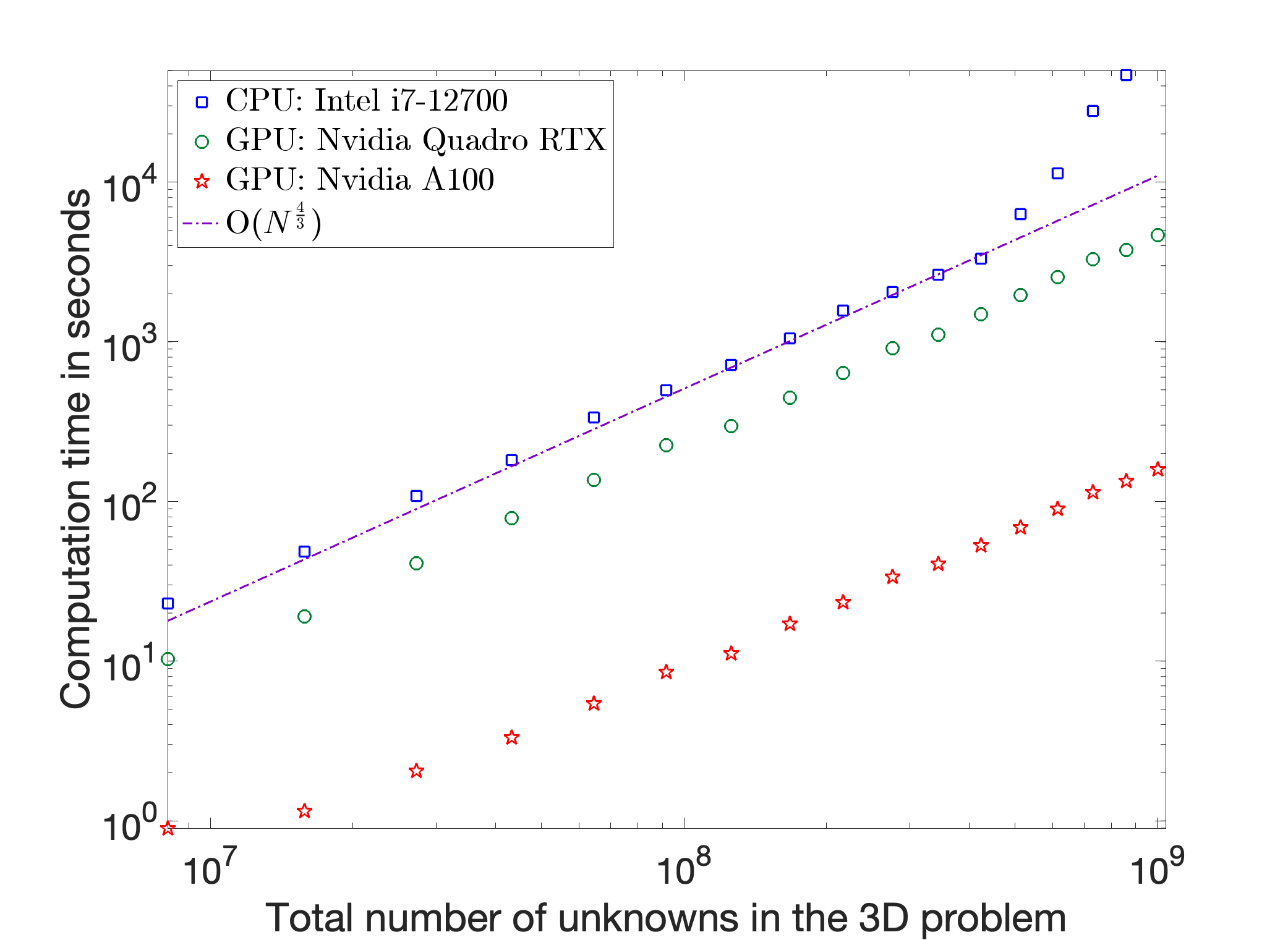}}
\caption{Online computation time of $Q^5$ spectral-element method solving a 3D Poisson equation two hundred times.
On the A100, it takes approximately only 0.8 second when solving one Poisson equation for the total number of DoFs being $1001^3$.}
	\label{fig1:poisson_cpu_gpu}
 \end{figure}

In this subsection, we list the online computational time comparison
for solving $\alpha u -\Delta u = f$ on $\Omega=[-1,1]^3$ with $\alpha=1$ and Neumann boundary conditions, by using the $Q^5$ spectral-element method. 
 To obtain a more accurate estimate of the online computational time, we count the online computation time for solving the Poisson equation   200 times.   The results of online computational time, depicted in Figure \ref{fig1:poisson_cpu_gpu} and Table \ref{tab2:poisson_cpu_gpu}, demonstrate a speed-up factor of at least 60 for sufficiently large problems  when comparing Nvidia A100 to Intel i7-12700. In particular, we observe that on the A100, solving a Poisson type equation \eqref{pde} with a total degree of freedoms (DoFs) equal to $1001^3$, takes approximately only 0.8 second.

For completeness, in Table \ref{tab7:poisson_cpu_gpu_offline}, we also include the offline preparation time which includes the time for generating arrays and loading arrays to GPU memory.

\begin{table}[htbp]
    \centering
    \begin{tabular}{|c|c|c c|c c|}
    \hline
     \multirow{2}{*}{Total DoFs}& Intel i7-12700 & \multicolumn{2}{c|}{NVIDIA Quadro} & \multicolumn{2}{c|}{NVIDIA A100}\\
     \cline{2-6} 
          & CPU time  & GPU time & speed-up & GPU time & speed-up\\
    \hline
          $201^3$ & $2.29$E$1$ & $1.03$E$1$ & $2.23$ & $9.00$E-$1$& $25.47$ \\
    \hline
          $251^3$ & $4.86$E$1$ & $1.91$E$1$ & $2.54$ & $1.15$E$0$ & $42.14$ \\
    \hline
          $301^3$ & $1.08$E$2$ & $4.10$E$1$ & $2.65$ & $2.05$E$0$ & $52.94$ \\
    \hline
          $351^3$ & $1.81$E$2$ & $7.82$E$1$ & $2.32$ & $3.31$E$0$ & $54.77$ \\
    \hline
          $401^3$ & $3.36$E$2$ & $1.36$E$2$ & $2.47$ & $5.41$E$0$ & $62.12$ \\
    \hline
          $451^3$ & $4.98$E$2$ & $2.25$E$2$ & $2.22$ & $8.52$E$0$ & $58.49$ \\
    \hline
          $501^3$ & $7.13$E$2$ & $2.96$E$2$ & $2.41$ & $1.11$E$1$ & $64.09$ \\
    \hline
          $551^3$ & $1.05$E$2$ & $4.46$E$2$ & $2.35$ & $1.71$E$1$ & $61.19$ \\
    \hline
          $601^3$ & $1.57$E$3$ & $6.40$E$2$ & $2.46$ & $2.35$E$1$ & $67.09$ \\
    \hline
          $651^3$ & $2.05$E$3$ & $9.09$E$2$ & $2.25$ & $3.37$E$1$ & $60.68$ \\
    \hline
          $701^3$ & $2.63$E$3$ & $1.11$E$3$ & $2.37$ & $4.07$E$1$ & $64.63$ \\
    \hline
          $751^3$ & $3.30$E$3$ & $1.48$E$3$ & $2.23$ & $5.31$E$1$ & $62.19$ \\
    \hline
          $801^3$ & $6.29$E$3$ & $1.97$E$3$ & $3.20$ & $6.86$E$1$ & $91.64$ \\
    \hline
          $851^3$ & $1.13$E$4$ & $2.55$E$3$ & $4.45$ & $8.97$E$1$ & $126.36$ \\
    \hline
          $901^3$ & $2.79$E$4$ & $3.27$E$3$ & $8.53$ & $1.14$E$2$ & $244.19$ \\
    \hline
          $951^3$ & $4.69$E$4$ & $3.77$E$3$ & $12.45$ & $1.34$E$2$ & $349.16$ \\
    \hline
    \end{tabular}
    \caption{\it Online computation time of solving a 3D Poisson equation two hundred times on three devices: the time unit is second, and the speed-up is GPU versus CPU. }
    \label{tab2:poisson_cpu_gpu}
\end{table}

\begin{table}[ht!]
    \centering
    \begin{tabular}{|c|c|c|c|c|c|c|}
    \hline
         Total DoFs  & $200^3$ & $250^3$ & $300^3$ & $350^3$ & $400^3$ & $450^3$\\
    \hline
         Quadro & $1.59$E$0$ & $1.60$E$0$ & $1.64$E$0$ & $1.72$E$0$ & $1.77$E$0$ & $1.85$E$0$\\
    \hline
         A100 & $3.01$E-$1$ & $3.19$E-$1$ & $3.41$E-$1$ & $3.85$E-$1$ & $4.36$E-$1$ & $4.93$E-$1$\\
    \hline
    \hline
         Total DoFs & $500^3$ & $550^3$ & $600^3$ & $650^3$ & $700^3$ & $750^3$\\
    \hline
         Quadro & $1.92$E$0$ & $2.01$E$0$ & $2.16$E$0$ & $2.42$E$0$ & $2.46$E$0$ & $2.49$E$0$ \\
    \hline
         A100 & $5.57$E-$1$ & $6.30$E-$1$ & $7.07$E-$1$ & $8.41$E-$1$ & $9.61$E-$1$ & $1.12$E$0$\\
    \hline
    \hline
         Total DoFs & $800^3$ & $850^3$ & $900^3$ & $950^3$ & $1000^3$ & $1050^3$\\
    \hline
         Quadro & $2.52$E$0$ & $2.82$E$0$ & $3.08$E$0$ & $3.40$E$0$ & $3.70$E$0$ & $4.37$E$0$\\
    \hline
         A100 & $1.24$E$0$ & $1.41$E$0$ & $1.60$E$0$ & $1.83$E$0$ & $2.04$E$0$ & $2.29$E$0$\\
    \hline
    \end{tabular}
    \caption{\it Offline preparation time in MATLAB on GPU for solving a 3D Poisson  equation: the time unit is second.}
    \label{tab7:poisson_cpu_gpu_offline}
\end{table}

\subsection{GPU acceleration for solving a Schr{\"o}dinger equation}
For a Problem with general variable coefficients,   the tensor product structure of the eigenvectors no longer holds. Then, an efficient method for solving such problems is to use a preconditioned conjugate gradient method with the inverse of Poisson type equation as   a preconditioner.
As an example,
we consider the following equation
\begin{equation}
\label{shrodinger}
   \alpha u -\Delta u + V(\bx)u = f,
\end{equation}
on $\Omega=[-16,16]^3$ with   $\alpha = 1$, 
\begin{equation}
\label{beta}
    V(\bx) = \beta\sin(\frac{\pi}{4}x)^2\sin(\frac{\pi}{4}y)^2\sin(\frac{\pi}{4}z)^2,\quad \beta>0,
\end{equation}
and an exact solution 
\begin{equation}
    u(\bx) = \cos(\frac{\pi}{16} x)\cos(\frac{\pi}{16} y)\cos(\frac{\pi}{16} z).
\end{equation}
The equation \eqref{shrodinger} is sometimes referred to as a Schr{\"o}dinger equation, which emerges in solving more complicated problems originated from the nonlinear Schr{\"o}dinger equation, e.g., the Gross-Pitaevskii equation \cite{chen2023convergence}.
The boundary conditions can be either periodic or homogeneous Neumann. 

Note that $0\le V(x)\le \beta$. 
We use $((\alpha+\frac 12\beta) I - \Delta)^{-1}$ as a preconditioner in the preconditioned conjugate gradient (PCG) method inverting the operator $\alpha I - \Delta +V(\mathbf x)$ with periodic boundary conditions in the $Q^5$ spectral-element method, where $((\alpha+\frac 12\beta) I - \Delta)^{-1}$ is implemented in the same way as in Table \ref{tab:3dcode}. {\bf We emphasize that eigenvectors can not be implemented by fast Fourier transform (FFT) for high order schemes with periodic boundary conditions, because the stiffness matrix $S$ for $
Q^k$ SEM is a circulant matrix only when $k=1$, i.e., FFT can be used to invert Laplacian only for second order accurate schemes.}

Obviously, the performance of such a simple method depends on the condition number of the operator $\alpha I - \Delta +V(\mathbf x)$, which is affected by the choice of $V(\mathbf x)$.  By choosing different $\beta$ in \eqref{beta}, the performance of PCG, e.g., the number of PCG iterations needed for the PCG iteration residue to reach round-off errors, would vary. We first list the performance of PCG for the $Q^5$ spectral-element method on different meshes for different $\beta$ in Table \ref{tab7c: schrodinger_iter}. We can observe that the performance only depends on $V(\bx)$ for a fine enough mesh.

\begin{table}[ht!]
    \centering
    \begin{tabular}{|c|c|c|c|c|c|c|c|c|c|c|c|}
    \hline
    Total&\multicolumn{10}{c|}{Number of PCG iterations}\\
    \cline{2-11}
        DoFs   & $\beta=1$  & $10$ & $100$& $200$& $400$ & $800$ & $1000$& $2000$& $4000$& $10000$\\
    \hline
          $250^3$ & $10$ & $35$ & $85$ & $112$& $149$ & $191$ & $214$ & $288$ & $388$ & $535$\\
    \hline
          $350^3$ & $10$ & $32$ & $81$ & $108$ & $143$ & $184$ & $202$ & $265$ & $348$ & $522$ \\
    \hline
          $450^3$ & $10$ & $30$ & $80$ & $105$ & $142$ & $181$ & $191$ & $252$ & $333$ & $467$ \\
    \hline
          $550^3$ & $10$ & $28$ & $76$ & $104$ & $129$ & $174$ & $183$ & $239$ & $321$ & $448$ \\
    \hline
          $650^3$ & $10$ & $27$ & $77$ & $100$ & $135$ & $169$ & $182$ & $248$ & $327$ & $465$ \\
    \hline
          $750^3$ & $10$ & $27$ & $76$ & $100$ & $130$ & $173$ & $192$ & $233$ & $320$ & $456$ \\
    \hline
          $850^3$ & $10$ & $27$ & $76$ & $100$ & $134$ & $173$ & $188$ & $233$ & $305$ & $430$ \\
    \hline
          $950^3$ & $10$ & $25$ & $72$ & $98$ & $128$ & $165$ & $180$ & $234$ & $305$ & $428$\\
    \hline
          $1000^3$ & $10$ & $25$ & $72$ & $101$ & $134$ & $172$ & $188$ & $243$ & $305$ & $433$\\
    \hline
    \end{tabular}
    \caption{\it Number of PCG iterations needed for PCG with $((\alpha+\frac{1}{2}\beta)I - \Delta)^{-1}$ as the preconditioner to converge for   solving a Schr{\"o}dinger equation by the $Q^5$ SEM on different meshes with different $\beta=1,10,100,200,\cdots, 10000$.}
    \label{tab7c: schrodinger_iter}
\end{table}

\begin{figure}[htbp]
 \subfigure[Comparison on three devices.]{\includegraphics[scale=0.18]{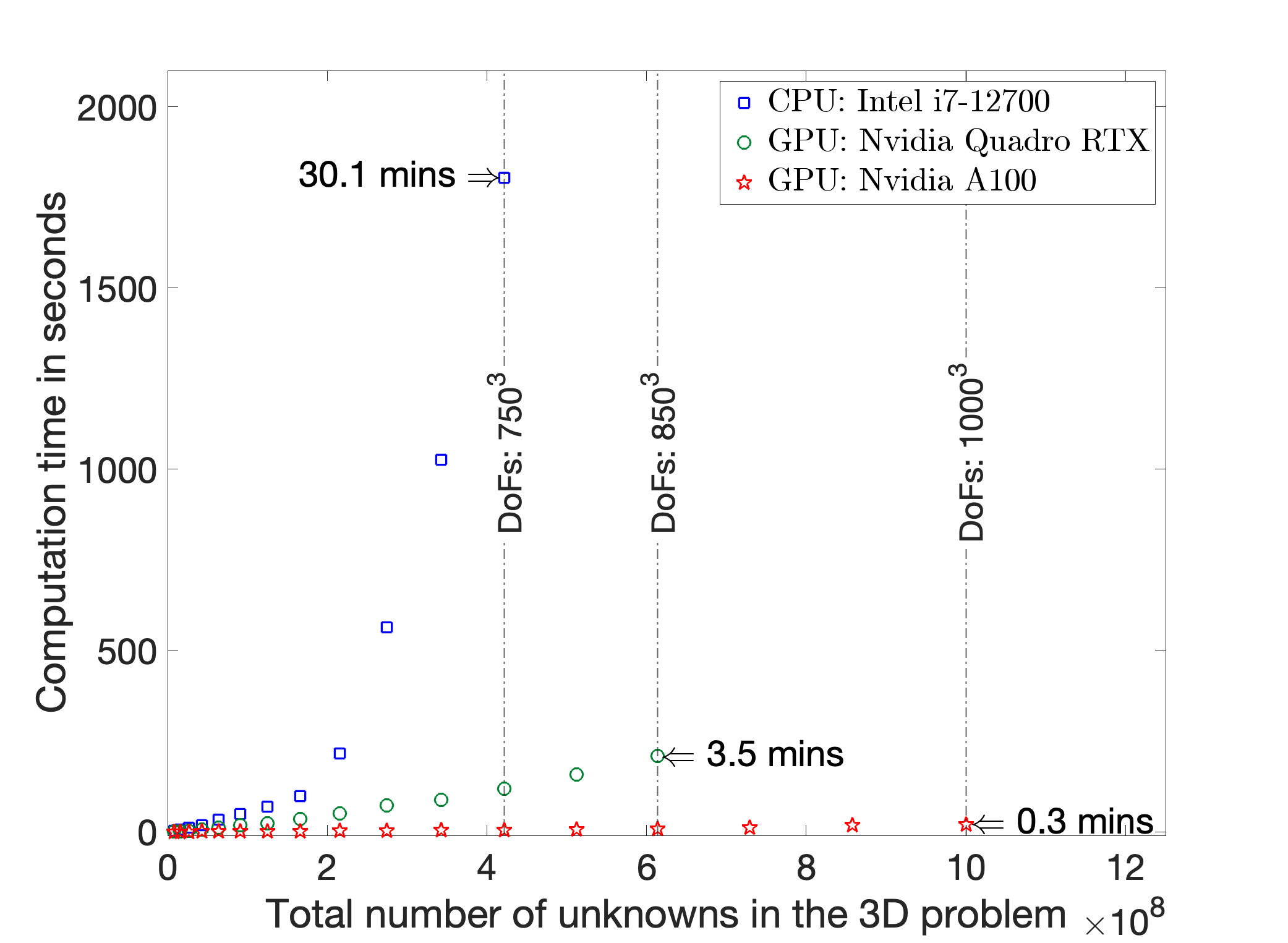}}
 \subfigure[Semilogx plot shows the complexity on all devices are $\mathcal O(N^{\frac43})$  for problems of proper sizes.]{\includegraphics[scale=0.18]{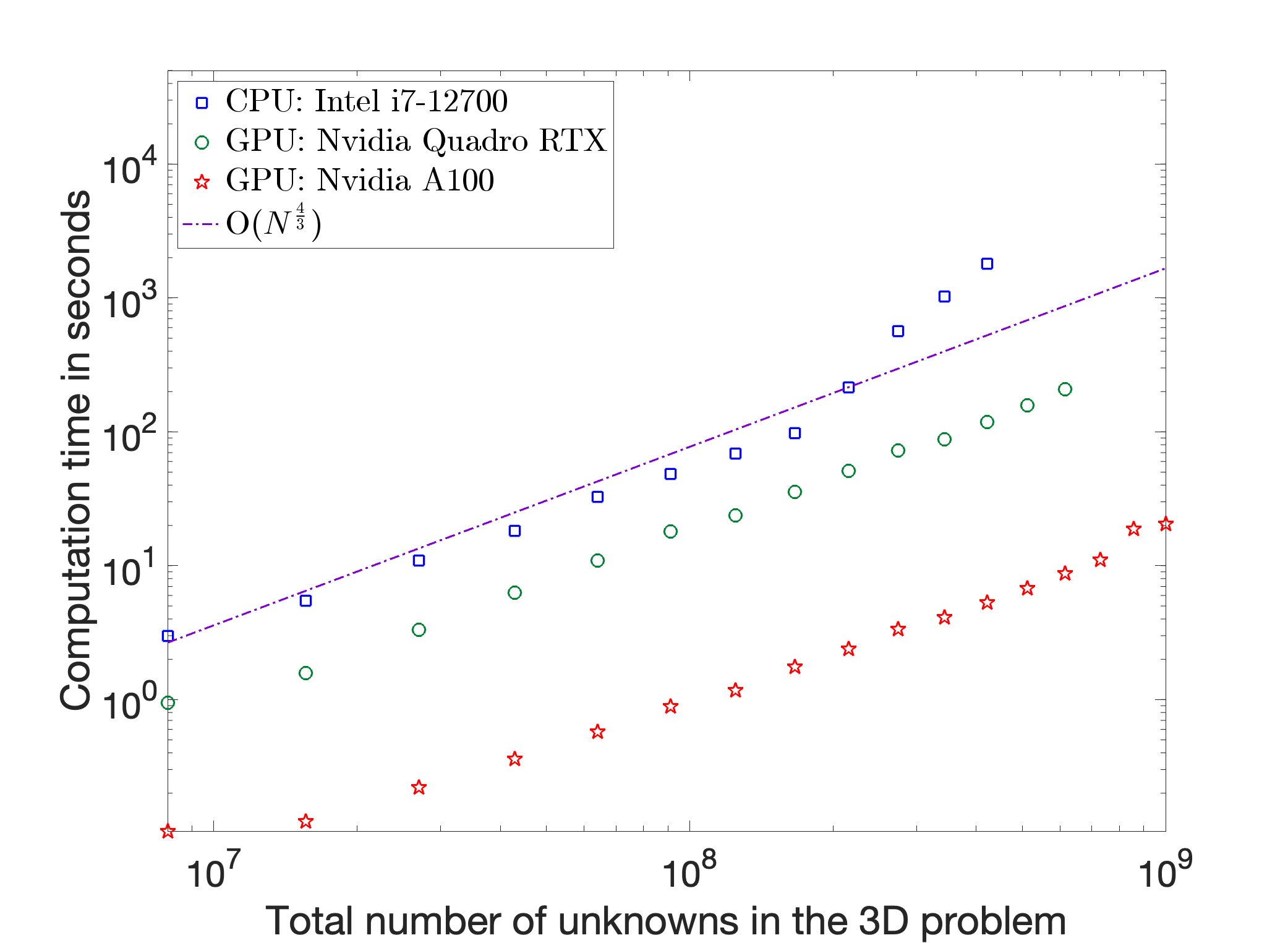}}
\caption{Online computational time of $Q^5$ SEM for a Schr{\"o}dinger equation with $\beta=1$ in \eqref{beta}, solved by PCG with $(I - \Delta)^{-1}$ as the preconditioner.}
	\label{fig1:schrodinger_cpu_gpu}
 \end{figure}

The online computational time of using PCG for the $Q^5$ spectral-element method solving one  Schr{\"o}dinger equation with $\beta=1$ in \eqref{beta} is listed in both Figure \ref{fig1:schrodinger_cpu_gpu} and Table \ref{tab5: schrodinger_CPU_GPU_beta1}. We can observe a satisfying speed-up. With 10 PCG iterations, it  costs about $20$ seconds on A100 for inverting a 3D Schr{\"o}dinger operator for a total number of DoFs as large as $1000^3$.

\begin{table}[ht!]
    \centering
    \begin{tabular}{|c|c|c|c|c|c|c|}
    \hline
         Total DoFs  & $200^3$ & $250^3$ & $300^3$ & $350^3$ & $400^3$ & $450^3$\\
    \hline
         Intel i7-12700 & $2.99$E$0$ & $5.50$E$0$ & $1.09$E$1$ & $1.81$E$1$ & $3.27$E$1$ & $4.83$E$1$\\
    \hline
         Nvidia Quadro & $9.48$E-$1$ & $1.58$E$0$ & $3.32$E$0$ & $6.30$E$0$ & $1.10$E$1$ & $1.80$E$1$\\
    \hline
         Nvidia A100    & $1.04$E-$1$ & $1.23$E-$1$ & $2.21$E-$1$ & $3.61$E-$1$ & $5.75$E-$1$ & $8.86$E-$1$\\
    \hline
    \hline
         Total DoFs & $500^3$ & $550^3$ & $600^3$ & $650^3$ & $700^3$ & $750^3$\\
    \hline
         Intel i7-12700 & $6.90$E$1$ & $9.83$E$1$ & $2.16$E$2$ & $5.63$E$2$ & $1.03$E$3$ & $1.80$E$3$\\
    \hline
         Nvidia Quadro & $2.38$E$1$ & $3.55$E$1$ & $5.11$E$1$ & $7.23$E$1$ & $8.78$E$1$ & $1.18$E$2$\\
    \hline
         Nvidia A100 & $1.18$E$0$ & $1.75$E$0$ & $2.38$E$0$ & $3.35$E$0$ & $4.10$E$0$ & $5.30$E$0$\\
    \hline
    \hline
         Total DoFs & $800^3$ & $850^3$ & $900^3$ & $950^3$ & $1000^3$ & \\
    \hline
      Nvidia   Quadro & $1.57$E$2$ & $2.09$E$2$ & - & - & - & \\
    \hline
      Nvidia   A100 & $6.79$E$0$ & $8.71$E$0$ & $1.11$E$1$ & $1.88$E$1$ & $2.04$E$1$ & \\
    \hline
    \end{tabular}
    \caption{\it The online computational time (unit is second) of using PCG for the $Q^5$ SEM  solving one  Schr{\"o}dinger equation with $\beta=1$ in \eqref{beta}.}
    \label{tab5: schrodinger_CPU_GPU_beta1}
\end{table}

\subsection{Robustness of the implementation for very high order elements}
\label{sec:robustness}

For very high order elements, it is important to have a robust procedure for finding the eigenvalue decomposition of the matrix $H$. We test the implementation in Remark \ref{remark:robustness} for the $Q^{20}$
spectral-element method solving the Schr{\"o}dinger equation.
 The error in Table \ref{tab8: Q20} and the online computational time in Table \ref{tab9: Q20} validate the robustness of the implementation. In other words, even for $Q^{20}$ element, the numerical computation of eigenvalue decomposition in Remark \ref{remark:robustness} is still accurate.

\begin{table}[ht!]
    \centering
    \begin{tabular}{|c|c|c|c|}
    \hline
     \multirow{2}{*}{Total DoFs}& \multicolumn{3}{c|}{$\ell^\infty$ error} \\ \cline{2-4}
        &$\beta = 1$ & $\beta = 10$ & $\beta = 100$\\ 
    \hline
         $500^3$ & $1.89$E-$13$ & $1.62$E-$13$ & $1.29$E-$13$ \\
    \hline
         $800^3$ &  $4.86$E-$13$ & $4.09$E-$13$ & $2.97$E-$13$\\
    \hline
         $1000^3$ &  $6.28$E-$13$ & $5.23$E-$13$ & $3.76$E-$13$ \\
    \hline
    \end{tabular}
    \caption{\it The $\ell^\infty$ error for $Q^{20}$ SEM solving one  Schr{\"o}dinger equation with different $\beta$ in \eqref{beta}. }
    \label{tab8: Q20}
\end{table}

\begin{table}[ht!]
    \centering
    \begin{tabular}{|c|c|c|c|c|c|c|}
    \hline
     \multirow{2}{*}{Total DoFs}& \multicolumn{2}{c|}{$\beta = 1$} & \multicolumn{2}{c|}{$\beta = 10$} & \multicolumn{2}{c|}{$\beta = 100$}\\ \cline{2-7}
         & Time & $\#$ PCG & Time & $\#$ PCG & Time & $\#$ PCG\\
    \hline
         $500^3$ & $1.23$E$0$ & 10 & $2.36$E$0$ & 23 & $6.47$E$0$ & 68\\
    \hline
         $800^3$ &  $7.13$E$0$ & 11 & $3.15$E$1$ &  56 & $7.84$E$1$ & 142\\
    \hline
         $1000^3$ &  $2.04$E$1$ & 10 & $4.40$E$1$ &  23 & $1.30$E$2$ & 70\\
    \hline
    \end{tabular}
    \caption{\it Online computational time in seconds for $Q^{20}$ SEM solving one  Schr{\"o}dinger equation with different $\beta$ in \eqref{beta}. }
    \label{tab9: Q20}
\end{table}

\subsection{Comparison with FFT on GPU}

It is also interesting to compare the implementation in Table \ref{tab:3dcode} with the performance of  fast Fourier transform (FFT) on GPU. 
In order to do so, we consider solving the  Poisson type  equation\eqref{pde} and the Schr{\"o}dinger equation \eqref{eqn: chschemesol}  with periodic boundary conditions using second order finite difference, or equivalently the $Q^1$ spectral-element method, for which the discrete Laplacian can be diagonalized by FFT, e.g., the eigenvector matrices $T^{-1}$ in \eqref{3D-Poisson} is the discrete Fourier transform matrix.
In other words, for $Q^1$ element with periodic boundary, the implementation in Table \ref{tab:3dcode}
can be replaced by the following implementation via FFT in MATLAB:
\begin{table}[htbp]
   \begin{lstlisting}
U = fftn(F);
U = U./Lambda3D;
U = real(ifftn(U));
    \end{lstlisting}
\caption{The FFT implementation of  a second order scheme for the Poisson equation  with periodic boundary conditions.}
\label{code:FFT}
\end{table}

\begin{figure}[ht!]
 \subfigure[Comparison on one device.]{\includegraphics[scale=0.18]{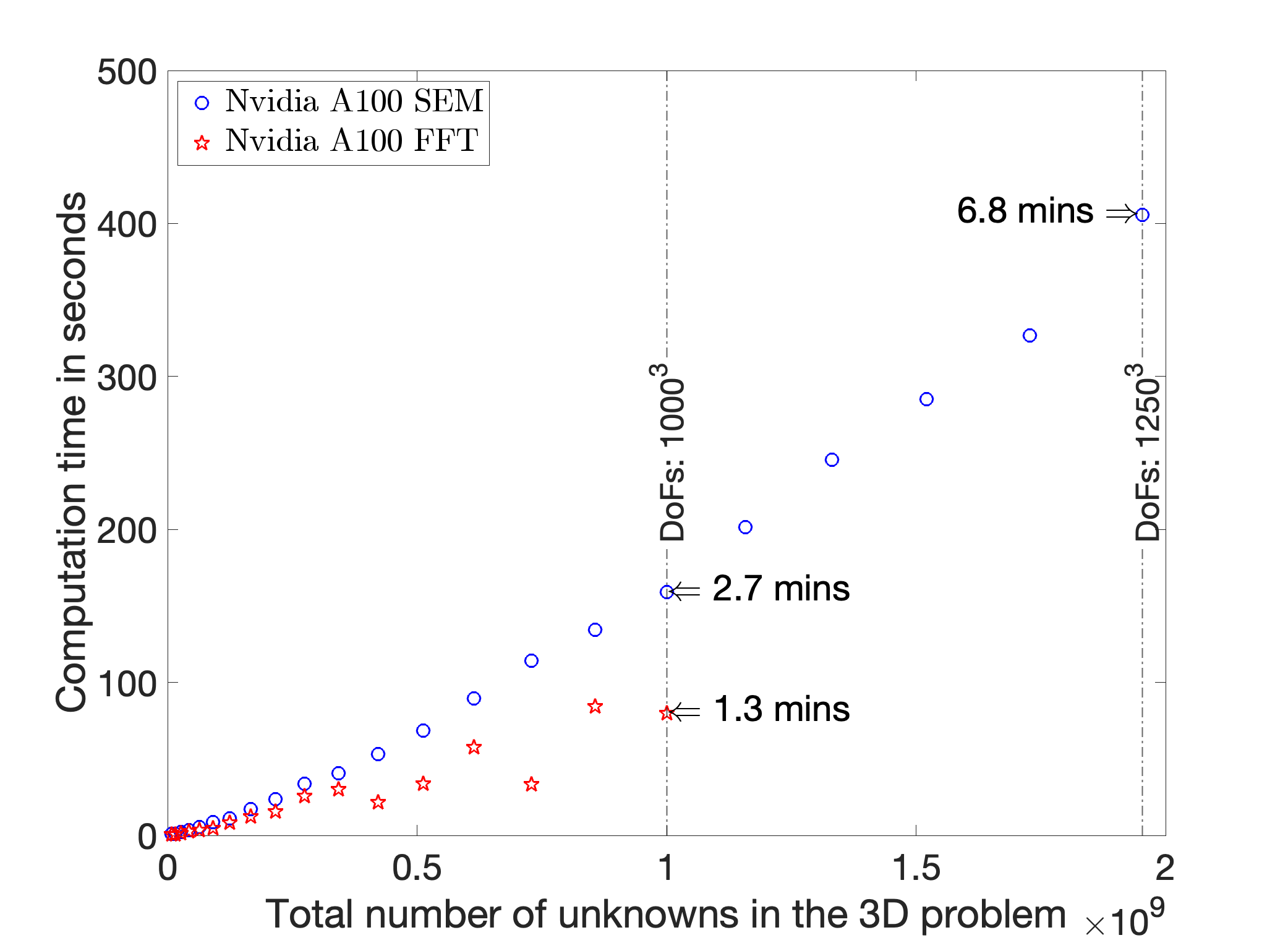}}
 \subfigure[Semilogx plot shows the complexity.]{\includegraphics[scale=0.18]{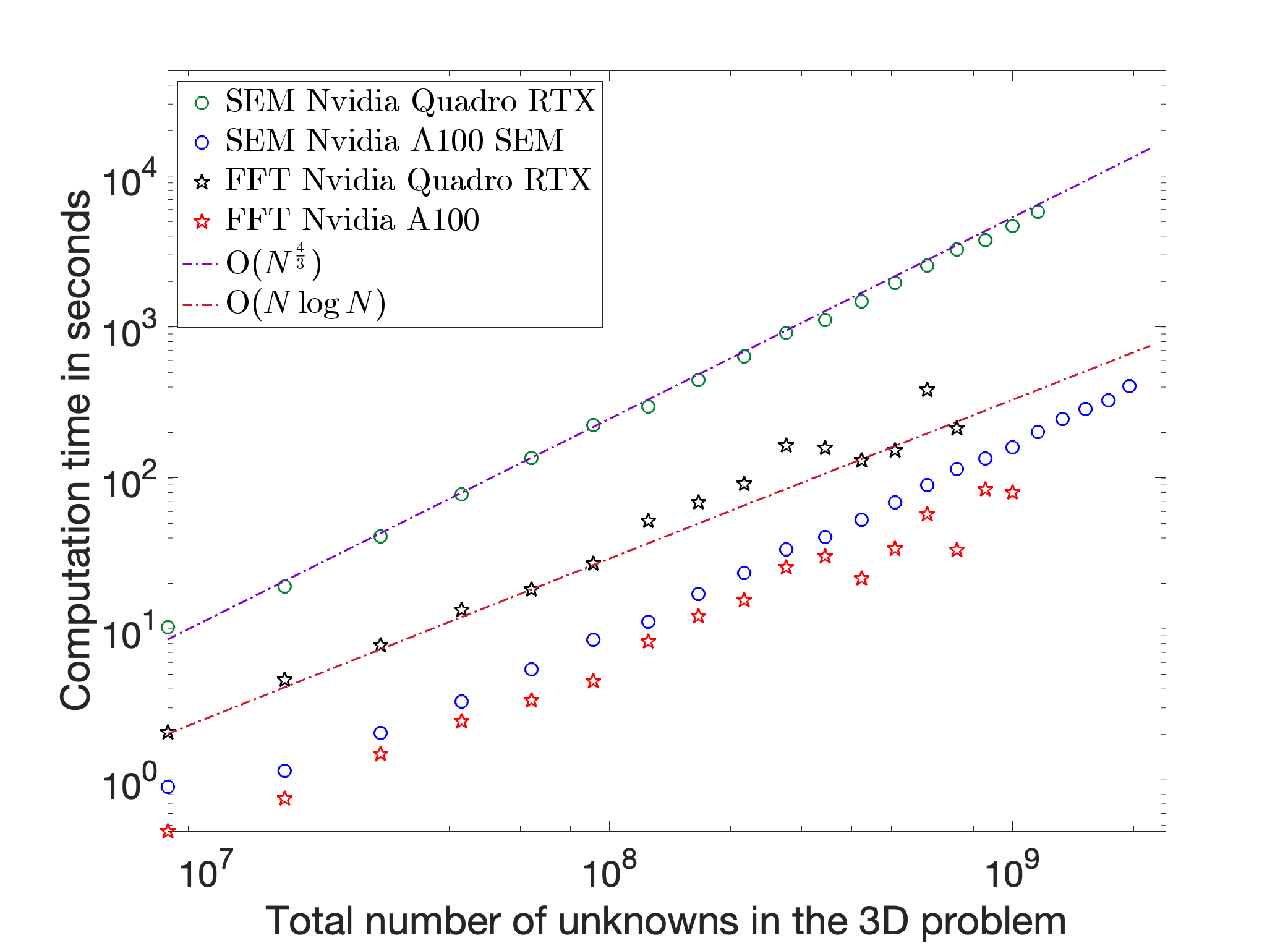}}
\caption{Comparison between $Q^5$ SEM implemented in Table \ref{tab:3dcode} and a second order scheme implemented by FFT in Table \ref{code:FFT}, for solving a Poisson equation 200 times. On A100, the FFT implementation cannot solve a problem of size $1050^3$ in MATLAB 2023, due to the larger memory cost of FFT. }
	\label{fig:fft1}
 \end{figure}

\begin{figure}[ht!]
 \subfigure[Comparison on one device.]{\includegraphics[scale=0.18]{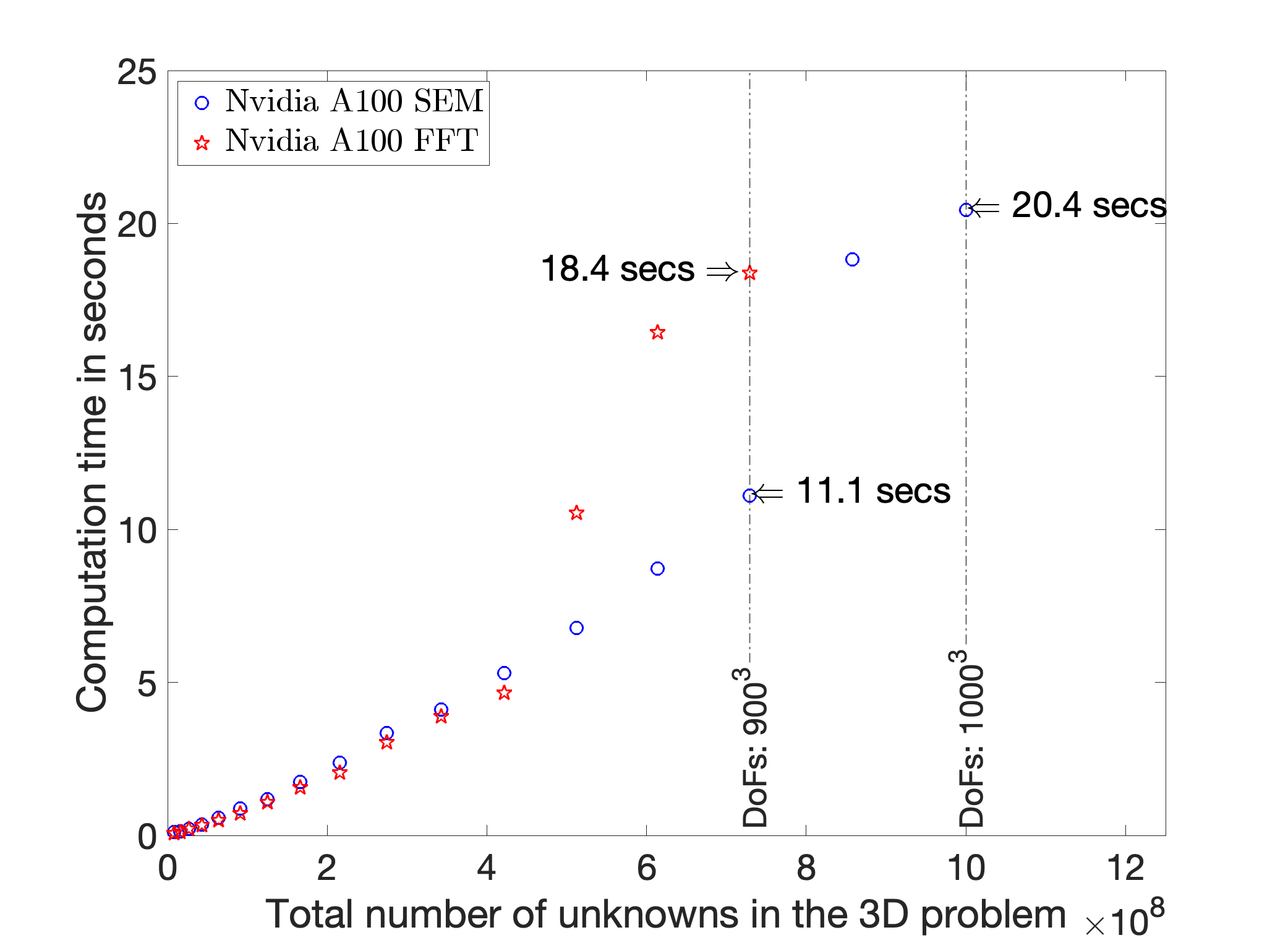}}
 \subfigure[Semilogx plot shows the complexity.]{\includegraphics[scale=0.18]{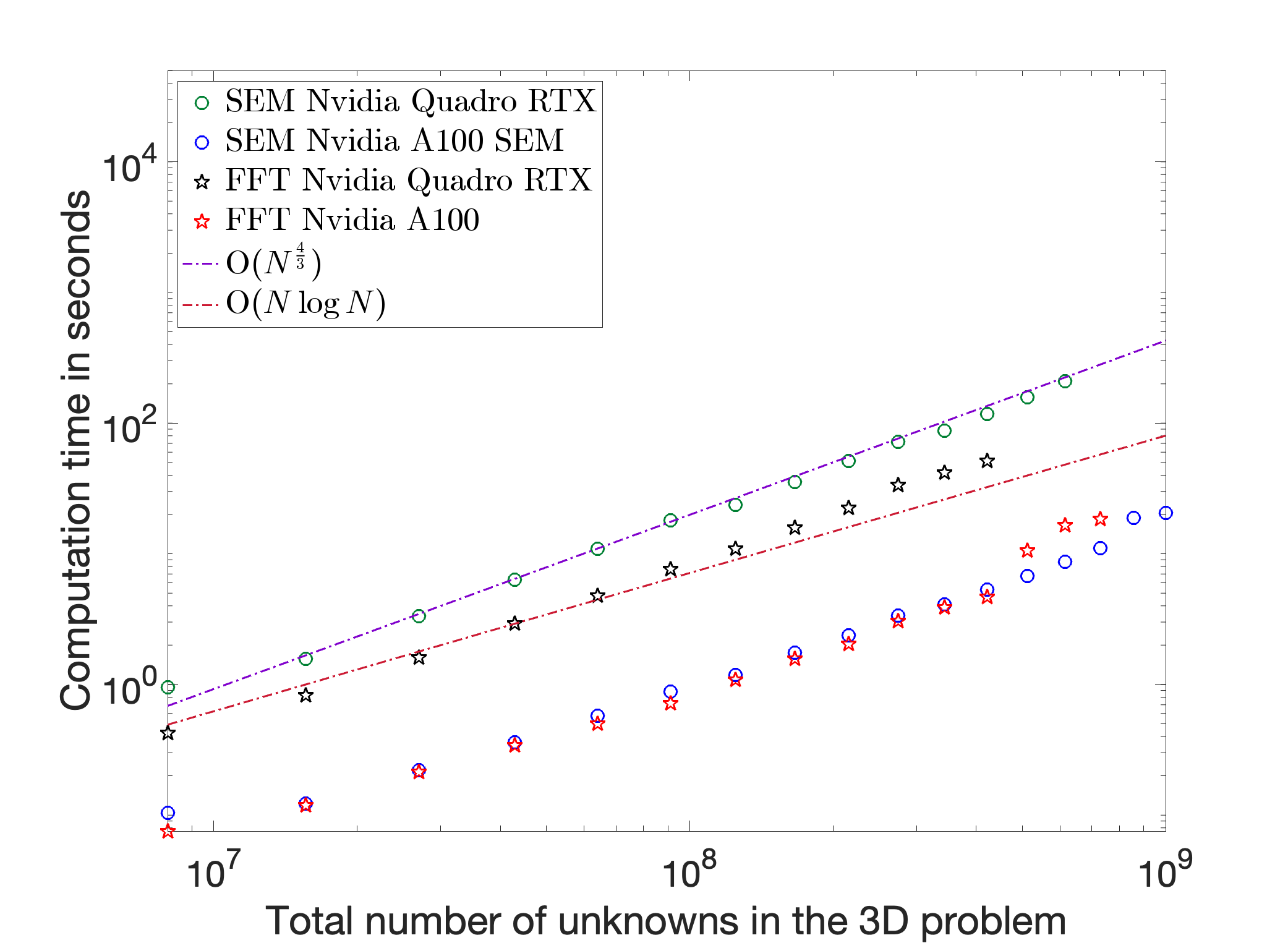}}
\caption{
Comparison between $Q^5$ SEM implemented in Table \ref{tab:3dcode} and a second order scheme implemented by FFT in Table \ref{code:FFT} for solving a Schr{\"o}dinger equation by PCG. On A100, the FFT implementation cannot solve a problem of size $950^3$ in MATLAB 2023, due to the larger memory cost of FFT.
}
	\label{fig:fft2}
 \end{figure}

 We will refer to such an implementation for a second order scheme as FFT in Figure \ref{fig:fft1} and  Figure \ref{fig:fft2}.
On the other hand, even if  the  Poisson equation has periodic boundary conditions, the matrices $T$ and $T^{-1}$ for high order elements cannot be implemented by FFT. We simply refer to the implementation in Table \ref{tab:3dcode} for  high order elements as SEM in  Figure \ref{fig:fft1} and  Figure \ref{fig:fft2}. 
The detailed comparison is listed in Figure \ref{fig:fft1} and  Figure \ref{fig:fft2}, as well as 
Table \ref{tab4: poisson_FFT_SEM} and Table \ref{tab6: schrodinger_FFT_SEM_beta1}. We can observe that FFT is faster as expected, on most meshes. However, the memory cost of performing FFT is more demanding, especially on finer meshes. For the Schr{\"o}dinger problem, the performance of FFT deteriorates on finest meshes.

\begin{table}[ht!]
    \centering
    \begin{tabular}{|c|c|c|c|c|c|c|}
    \hline
         Total DoFs & $200^3$ & $250^3$ & $300^3$ & $350^3$ & $400^3$ & $450^3$\\
    \hline
         Quadro (SEM) & $1.03$E$1$ & $1.91$E$1$ & $4.10$E$1$ & $7.82$E$1$ & $1.36$E$2$ &
         $2.25$E$2$\\
    \hline
         Quadro (FFT) & $2.06$E$0$ & $4.61$E$0$ & $7.82$E$0$ & $1.33$E$1$ & $1.83$E$1$ &
         $2.70$E$1$\\
    \hline
         A100 (SEM) & $9.00$E-$1$ & $1.15$E$0$ & $2.05$E$0$ & $3.31$E$0$ & $5.41$E$0$ &
         $8.52$E$0$\\
    \hline
         A100 (FFT) & $4.58$E-$1$ & $7.50$E-$1$ & $1.49$E$0$ & $2.44$E$0$ & $3.37$E$0$ &
         $4.52$E$0$\\
    \hline
    \hline
         Total DoFs & $500^3$ & $550^3$ & $600^3$ & $650^3$ & $700^3$ & $750^3$\\
    \hline
         Quadro (SEM) & $2.96$E$2$ & $4.46$E$2$ & $6.40$E$2$ & $9.09$E$2$ & $1.11$E$3$ &
         $1.48$E$3$\\
    \hline
         Quadro (FFT) & $5.19$E$1$ & $6.89$E$1$ & $9.12$E$1$ & $1.64$E$2$ & $1.57$E$2$ &
         $1.31$E$2$\\
    \hline
         A100 (SEM) & $1.11$E$1$ & $1.71$E$1$ & $2.35$E$1$ & $3.37$E$1$ & $4.07$E$1$ &
         $5.31$E$1$\\
    \hline
         A100 (FFT) & $8.22$E$0$ & $1.21$E$1$ & $1.56$E$1$ & $2.55$E$1$ & $3.02$E$1$ & $2.17$E$1$\\
    \hline
    \hline
         Total DoFs & $800^3$ & $850^3$ & $900^3$ & $950^3$ & $1000^3$ & $1050^3$\\
    \hline
         Quadro (SEM) & $1.97$E$3$ & $2.55$E$3$ & $3.27$E$3$ & $3.77$E$3$ & $4.656$E$3$ & $5.79$E$3$\\
    \hline
         Quadro (FFT) & $1.52$E$2$ & $3.82$E$2$ & $2.13$E$2$ & - & - & -\\
    \hline
         A100 (SEM) & $6.86$E$1$ & $8.97$E$1$ & $1.14$E$2$ & $1.34$E$2$ & $1.59$E$2$ &
         $2.01$E$2$\\
    \hline
         A100 (FFT) & $3.38$E$1$ & $5.76$E$1$ & $3.32$E$1$ & $8.43$E$1$ & $7.97$E$1$ & -\\
    \hline
    \hline
         Total DoFs & $1100^3$ & $1150^3$ & $1200^3$ & $1250^3$ & $1300^3$ & $1350^3$\\
    \hline
         A100 (SEM) & $2.46$E$2$ & $2.85$E$2$ & $3.27$E$2$ & $4.06$E$2$ & - & -\\
    \hline
    \end{tabular}
    \caption{\it Online computational time comparison between $Q^5$ SEM implemented in Table \ref{tab:3dcode} and a second order scheme implemented by FFT in Table \ref{code:FFT}, for solving a Poisson equation 200 times. On A100, the FFT implementation cannot solve a problem of size $1050^3$ in MATLAB 2023, due to the larger memory cost of FFT.  Unit is in seconds.}
    \label{tab4: poisson_FFT_SEM}
\end{table}
\begin{table}[ht!]
    \centering
    \begin{tabular}{|c|c|c|c|c|c|c|}
    \hline
         Total DoFs  & $200^3$ & $250^3$ & $300^3$ & $350^3$ & $400^3$ & $450^3$\\
    \hline
         Quadro (SEM) & $9.48$E-$1$ & $1.58$E$0$ & $3.32$E$0$ & $6.30$E$0$ & $1.10$E$1$ & $1.80$E$1$\\
    \hline
         Quadro (FFT) & $4.25$E-$1$ & $8.22$E-$1$ & $1.60$E$0$ & $2.91$E$0$ & $4.76$E$0$ & $7.60$E$0$\\
    \hline
         A100 (SEM) & $1.04$E-$1$ & $1.23$E-$1$ & $2.21$E-$1$ & $3.61$E-$1$ & $5.75$E-$1$ & $8.86$E-$1$\\
    \hline
         A100 (FFT) & $7.55$E-$2$ & $1.19$E-$1$ & $2.13$E-$1$ & $3.41$E-$1$ & $5.02$E-$1$ & $7.17$E-$1$\\
    \hline
    \hline
         Total DoFs & $500^3$ & $550^3$ & $600^3$ & $650^3$ & $700^3$ & $750^3$\\
    \hline
         Quadro (SEM) & $2.38$E$1$ & $3.55$E$1$ & $5.11$E$1$ & $7.23$E$1$ & $8.78$E$1$ & $1.18$E$2$\\
    \hline
         Quadro (FFT) & $1.09$E$1$ & $1.58$E$1$ & $2.24$E$1$ & $3.36$E$1$ & $4.16$E$1$ & $5.14$E$1$\\
    \hline
         A100 (SEM) & $1.18$E$0$ & $1.75$E$0$ & $2.38$E$0$ & $3.35$E$0$ & $4.10$E$0$ & $5.30$E$0$\\
    \hline
         A100 (FFT) & $1.09$E$0$ & $1.57$E$0$ & $2.05$E$0$ & $3.04$E$0$ & $3.89$E$0$ & $4.66$E$0$\\
    \hline
    \hline
         Total DoFs & $800^3$ & $850^3$ & $900^3$ & $950^3$ & $1000^3$ & $1050^3$\\
    \hline
         Quadro (SEM) & $1.57$E$2$ & $2.09$E$2$ & - & - & - & -\\
    \hline
         Quadro (FFT) & - & - & - & - & - & -\\
    \hline
         A100 (SEM) & $6.79$E$0$ & $8.71$E$0$ & $1.11$E$1$ & $1.88$E$1$ & $2.04$E$1$ & -\\
    \hline
         A100 (FFT) & $1.05$E$1$ & $1.64$E$1$ & $1.84$E$1$ & - & - & -\\
    \hline
    \end{tabular}
    \caption{\it Online computational time comparison between $Q^5$ SEM implemented in Table \ref{tab:3dcode} and a second order scheme implemented by FFT in Table \ref{code:FFT} for solving a Schr{\"o}dinger equation by PCG. On A100, the FFT implementation cannot solve a problem of size $950^3$ in MATLAB 2023, due to the larger memory cost of complex numbers. Unit is in seconds.}
    \label{tab6: schrodinger_FFT_SEM_beta1}
\end{table}
\newpage
\subsection{A Cahn--Hilliard equation}
We consider solving the Cahn--Hilliard equation \cite{cahn1958free}, which is not only a fourth-order equation in space, but also incorporates a time derivative. Consider a domain $\Omega=[-1,1]^3$ with its boundary denoted as $\partial \Omega$. Within this domain, the Cahn--Hilliard equation with simple boundary conditions is given by
\begin{equation}\label{CH}
    \begin{cases}
        \phi_t = m\Delta\left(-\epsilon \Delta \phi + \frac{1}{\epsilon}F'(\phi)\right)\quad\text{in}\,\,\Omega,\\
        \partial_{\bm{n}}\phi =0,\quad \partial_{\bm{n}}\Delta \phi = 0\quad\text{on}\,\,\partial\Omega,
    \end{cases}
\end{equation}
where $\phi$ is a phase function with a thin, smooth transitional layer, whose thickness is proportional to the parameter $\epsilon$, $m$ is the mobility constant, and $F(\phi)=\frac{1}{4}(\phi^2-1)^2$ is a double-well form function.

Due to the simplicity of the boundary conditions, 
we can avoid solving a fourth-order equation directly by reformulating \eqref{CH} as a system of second-order equations after introducing the chemical potential $\mu$, which can be expressed as the variational derivative of the energy functional:
 \begin{equation}
 	E(\phi)=\int_{\Omega} \frac{\epsilon}{2}|\nabla \phi|^2+\frac{1}{\epsilon}F(\phi)d\bm{x}.
 \end{equation}
Then, the system can be derived as
\begin{equation}\label{CH-sys}
    \begin{cases}
        \phi_t - m\Delta \mu = 0\quad\text{in}\,\,\Omega,\\
        \mu = -\epsilon \Delta \phi + \frac{1}{\epsilon}F'(\phi)\quad\text{in}\,\,\Omega,\\
        \partial_{\bm{n}}\phi =0,\quad \partial_{\bm{n}}\mu = 0\quad\text{on}\,\,\partial\Omega.
    \end{cases}
\end{equation}
 For the space discretization, we use $Q^5$ spectral-element method. For time discretization, we implement the second order backward differentiation formula (BDF-$2$) to the system \eqref{CH-sys}:
\begin{equation}\label{BDF2}
    \begin{cases}
        \frac{a\phi_{n+1} - \hat{\phi}_{n}}{\delta t} - m\Delta \mu_{n+1} = 0,\\
        \mu_{n+1} = -\epsilon\Delta \phi_{n+1} + \frac{1}{\epsilon}F'(\bar{\phi}_{n}),
    \end{cases}
\end{equation}
where $a=\frac{3}{2}$, $\hat{\phi}_n = 2\phi_n - \frac{1}{2}\phi_{n-1}$, and $\bar{\phi}_n = 2\phi_n - \phi_{n-1}$.
To solve this linear system, we can write it as
\begin{equation}
    \begin{bmatrix}
        \alpha I & -m\delta t\Delta\\
        \epsilon\Delta & I
    \end{bmatrix}
    \begin{bmatrix}
        \phi\\
        \mu
    \end{bmatrix}
    =
    \begin{bmatrix}
        f_1\\
        f_2
    \end{bmatrix},
\end{equation}
and its solution is given by
\begin{equation}
    \begin{bmatrix}
        \phi\\
        \mu
    \end{bmatrix}
    =
    \begin{bmatrix}
         \mathcal{D} I & \mathcal{D}(m\delta t\Delta)\\
        \mathcal{D} (-\epsilon\Delta) & \alpha \mathcal{D} I
    \end{bmatrix}
    \begin{bmatrix}
        f_1\\
        f_2
    \end{bmatrix}
    =
    \begin{bmatrix}
        \mathcal{D}(f_1 + m\delta t\Delta f_2)\\
        \mathcal{D}(-\epsilon\Delta f_1 + \alpha f_2)
    \end{bmatrix},
    \label{eqn: chschemesol}
\end{equation}
where $\mathcal{D} = (\alpha I + m\delta t\epsilon\Delta^2)^{-1}$. 

Notice that $\mu$ and $\phi$ are already decoupled in \eqref{eqn: chschemesol}. Thus for implementing the scheme \eqref{BDF2}, we only need to compute $\phi$ without computing $\mu$:
\begin{equation}
    \label{chscheme-2}
    \phi_{n+1} = \mathcal{D}  \hat{\phi}_n+m\delta t \frac{1}{\epsilon}\mathcal{D} \Delta F'(\bar{\phi}_{n}),
\end{equation}
where
both $\mathcal{D} = (\alpha I + m\delta t\epsilon\Delta^2)^{-1}$ and $\mathcal{D}\Delta = (\alpha I + m\delta t\epsilon\Delta^2)^{-1}\Delta$ can be implemented in the same way as shown in Table \ref{tab:3dcode}.

Since $(\alpha I + m\delta t\epsilon\Delta^2)^{-1}$ and $(\alpha I + m\delta t\epsilon\Delta^2)^{-1}\Delta$ share the same eigenvectors, the implementation of \eqref{chscheme-2} costs slightly less than solving the  Poisson type equation twice. In Table \ref{tab2:poisson_cpu_gpu},
we observe that,  the average online computational time of inverting Laplacian once is approximately $0.8$ second for the number DoFs being $1001^3$. For the same mesh and same DoFs, each time step \eqref{chscheme-2} of solving the Cahn--Hilliard equation costs about $1.27$ seconds in Table \ref{tab5: Cahn-Hilliard_GPU}.
\subsubsection{Accuracy test}
\begin{figure}[ht!]
\centering
    \begin{tikzpicture}
        \begin{loglogaxis}[
            scaled y ticks=false,  
            axis line style= semithick,
            height=6cm, width=6cm,
            xlabel={$\delta t$},
            xmin=5e-5 , xmax=1e-1,
            minor x tick num=5,
            ylabel={$Error$},
            ylabel style={rotate=0},
            ymin=1e-13,ymax=1e-3,
            minor y tick num=5,
            ytick={1e-12,1e-10,1e-8,1e-6,1e-4,1e-2,1e0},
            legend style={draw=none,legend cell align=left},
            legend pos=south east,
            ]
            \addplot[smooth, black] 
            table[x=dt, y expr=\thisrow{dt2}]{2nd_accurancy.data};
            \addlegendentry{slope 2 line};
            \addplot[smooth, black,mark=*, mark options={scale=0.6,fill=green!100}] 
            table[x=dt, y=err_phi_2nd]{2nd_accurancy.data};
            \addlegendentry{$\ell^2$ relative error};
        \end{loglogaxis}
    \end{tikzpicture}
    \caption{\it The $\ell^2$  {relative} error of BDF2 scheme \eqref{chscheme-2} for the Cahn--Hilliard equation.}
    \label{fig1:CH_acc}
\end{figure}
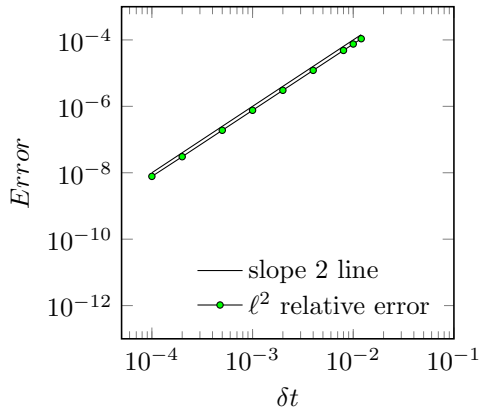

We first use a manufactured analytical solution of the Cahn--Hilliard equation to validate the convergence rate of the BDF-2 scheme \eqref{BDF2}. This solution is in the domain $\Omega=[-1,1]^3$ with $\epsilon = 0.2$, $m=0.01$: 
\begin{equation}
    \phi^{*}(\bx) =\cos(\pi x)\cos(\pi y)\cos(\pi z)\exp(t),
\end{equation}
and the corresponding forcing term can be obtained from the equation \eqref{CH-sys}. We fix the number of basis function as $N_x = N_y = N_z = 51$ in $Q^5$ SEM so that the spatial error is negligible compared with the time discretization error. Figure \ref{fig1:CH_acc} shows that the scheme \eqref{BDF2} achieves the expected second order time accuracy.
\subsubsection{Coalescence of two drops}\label{sec: cahn-hliiard}
We now study the coalescence of two droplets, as described by the Cahn--Hilliard equation, within the computational domain $\Omega=[-1,1]^3$.  Drawing from parameter settings in \cite{ChenShen2012}, we select $\epsilon=0.02$, the mobility constant $m=0.02$, and the time step size $\delta t = 0.001$ with an end time $T = 10$. For stable computation, we use the same simple stabilization method and stabilization parameter as in \cite{ChenShen2012}.
Initially, at time $t=0$, the domain is occupied by two neighboring spherical regions of the first material, while the second material fills the remaining space. As time progresses under the Cahn--Hilliard dynamics, these two spherical regions coalesce to form a singular droplet. More specifically, the initial condition for the phase function is given by
\begin{equation}
    \phi_{0}(\bx) = 1 - \tanh\frac{\vert \bx - \bx_1\vert - R}{\sqrt{2}\epsilon}  - \tanh\frac{\vert \bx - \bx_2\vert - R}{\sqrt{2}\epsilon},
\end{equation}
where $\bx_1 = (x_1,y_1,z_1) = (0,0,0.37)$ and $\bx_2 = (x_2,y_2,z_2) = (0,0,-0.37)$ are the centers of the initial spherical regions of the first material, and $R=0.35$ is the radius of these spheres.
\begin{figure}[ht!]
    \centering
    \includegraphics[width=0.85\linewidth]{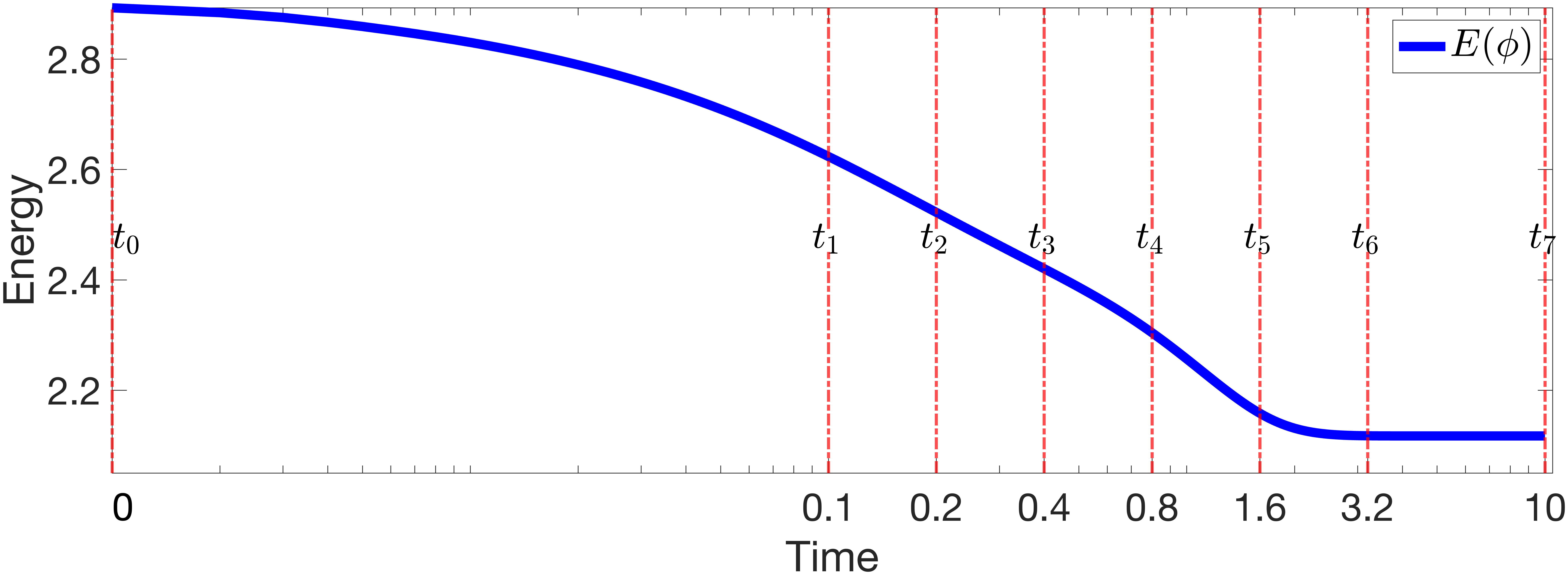}
    \caption{\it Semilogx plot shows the temporal evolution of energy $E(\phi)$.}
    \label{fig2:CH_energy}
\end{figure}

\begin{figure}[ht!]
    \begin{center}
	   \begin{tabular}{cccc}
        $t_{0}=0$  & $t_{1}=0.1$  & $t_{2}=0.2$ & $t_{3}=0.4$\\
        \includegraphics[width=0.21\linewidth]{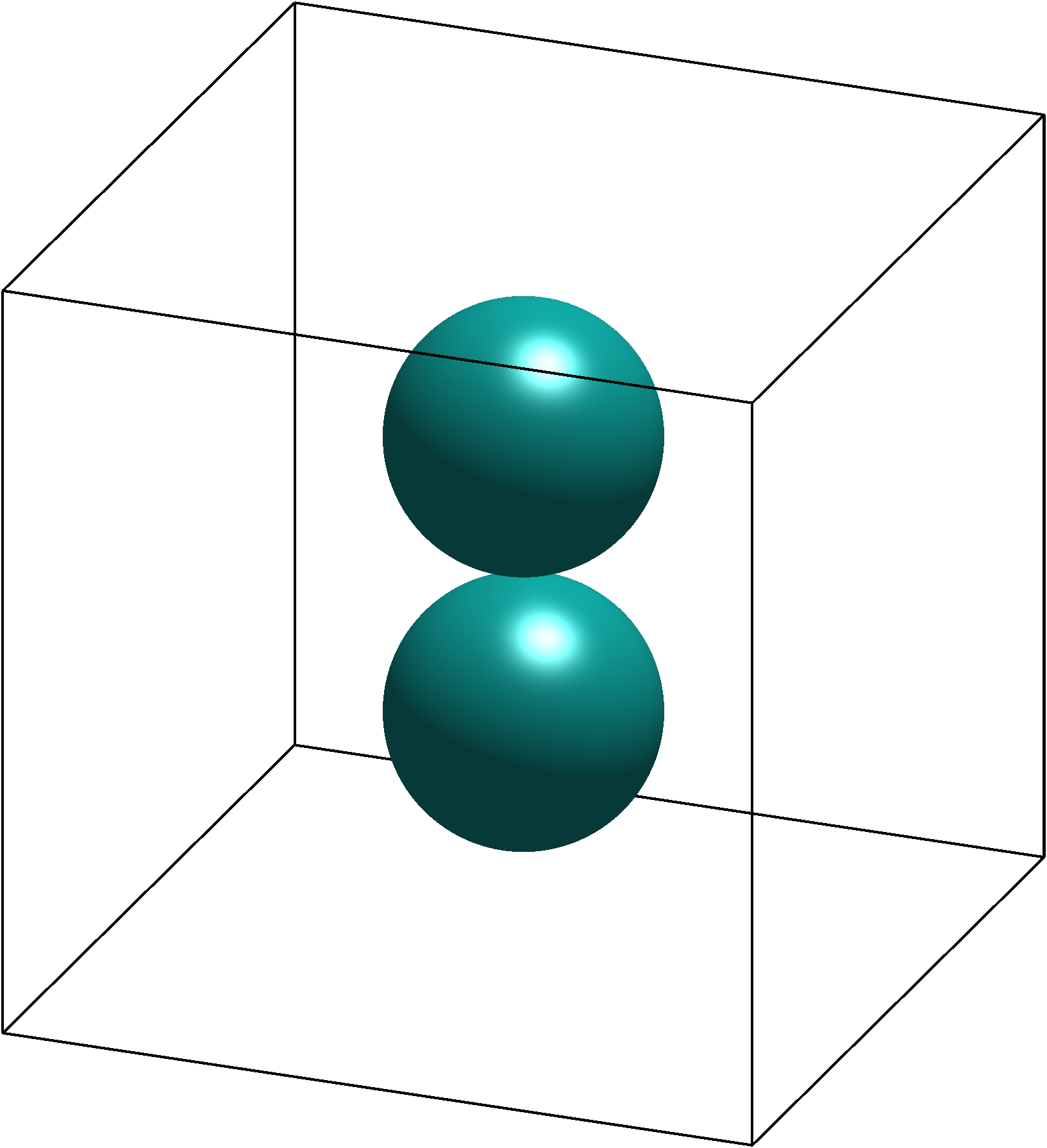}&
        \includegraphics[width=0.21\linewidth]{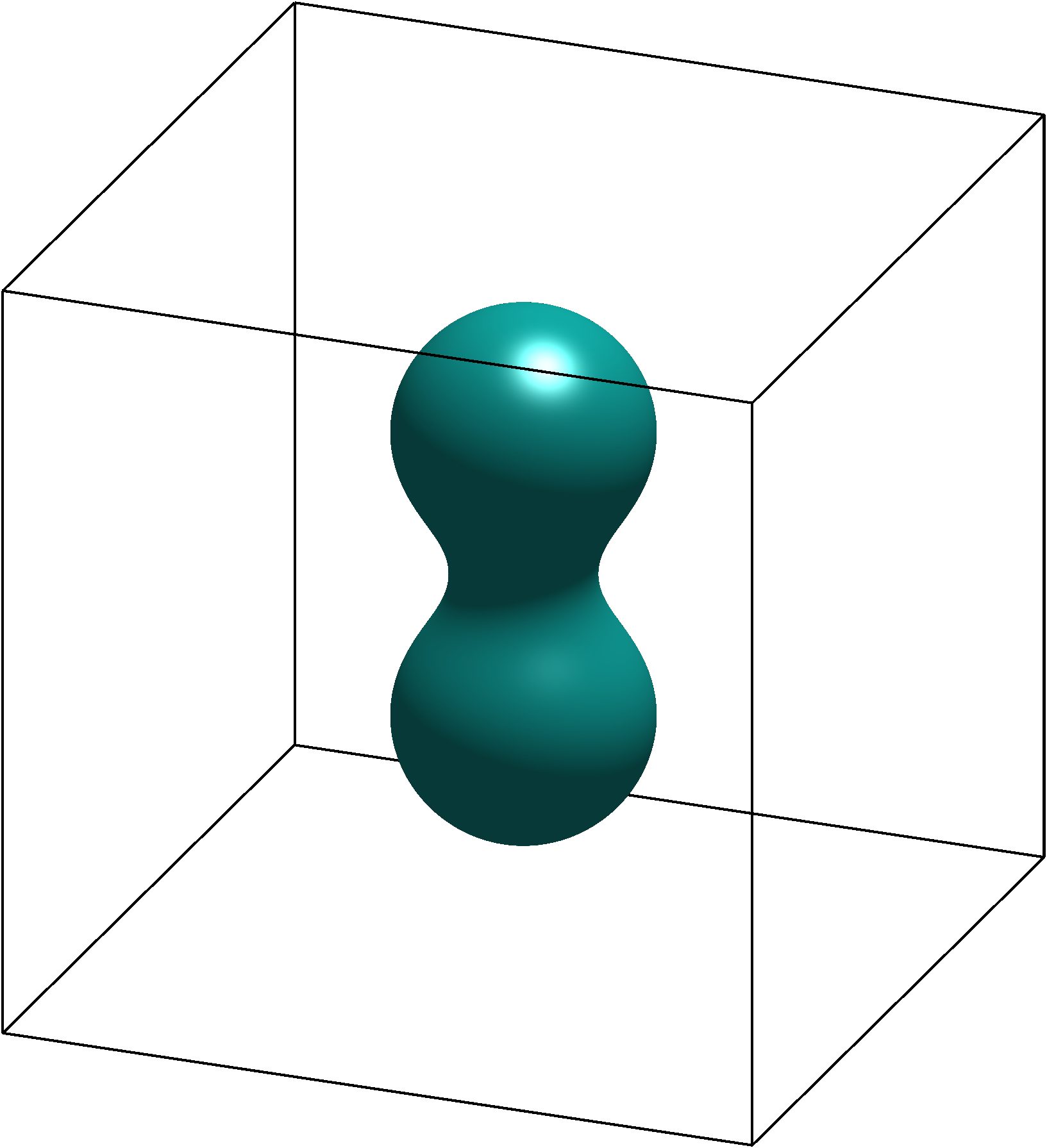}&
        \includegraphics[width=0.21\linewidth]{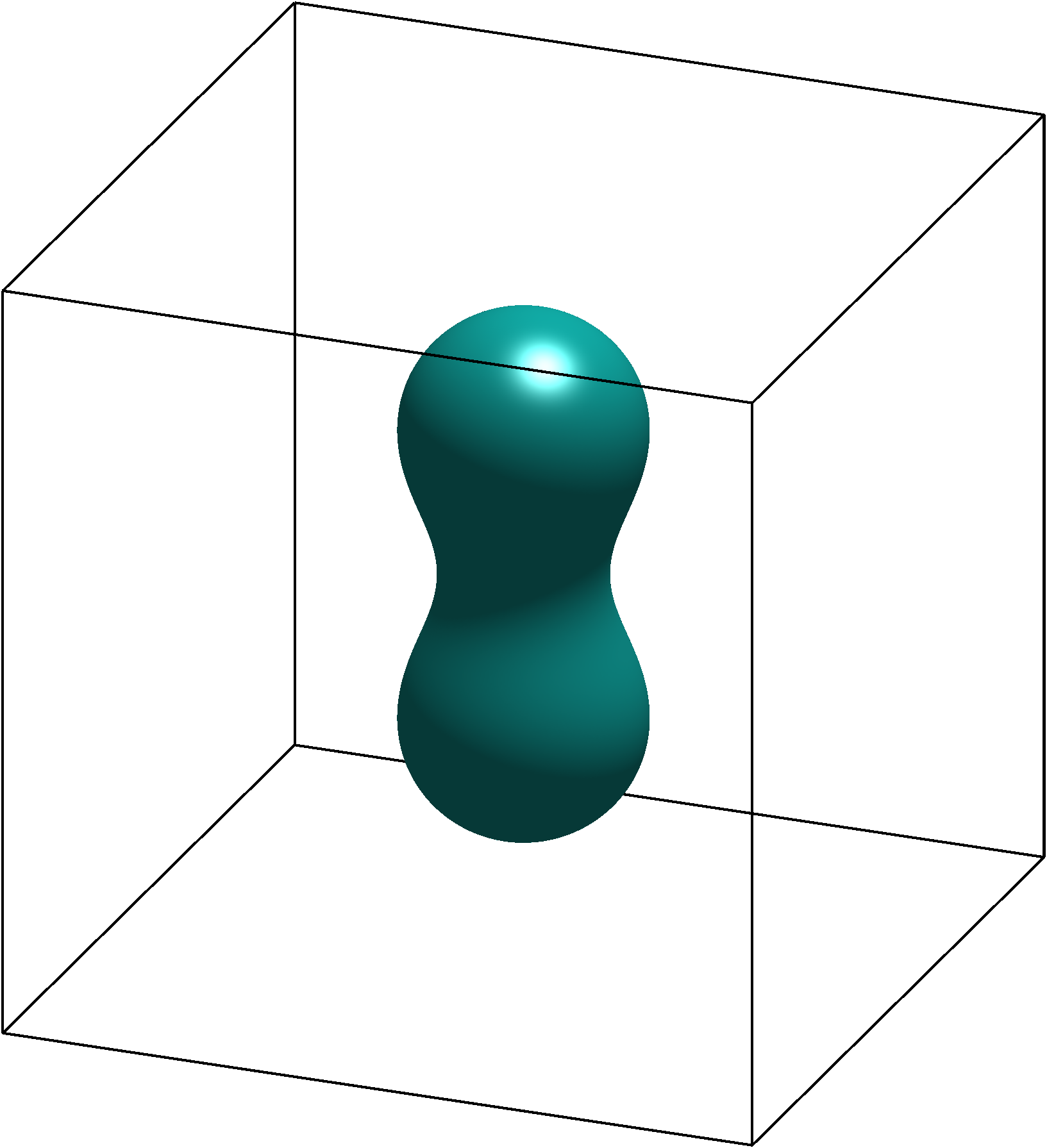}&
        \includegraphics[width=0.21\linewidth]{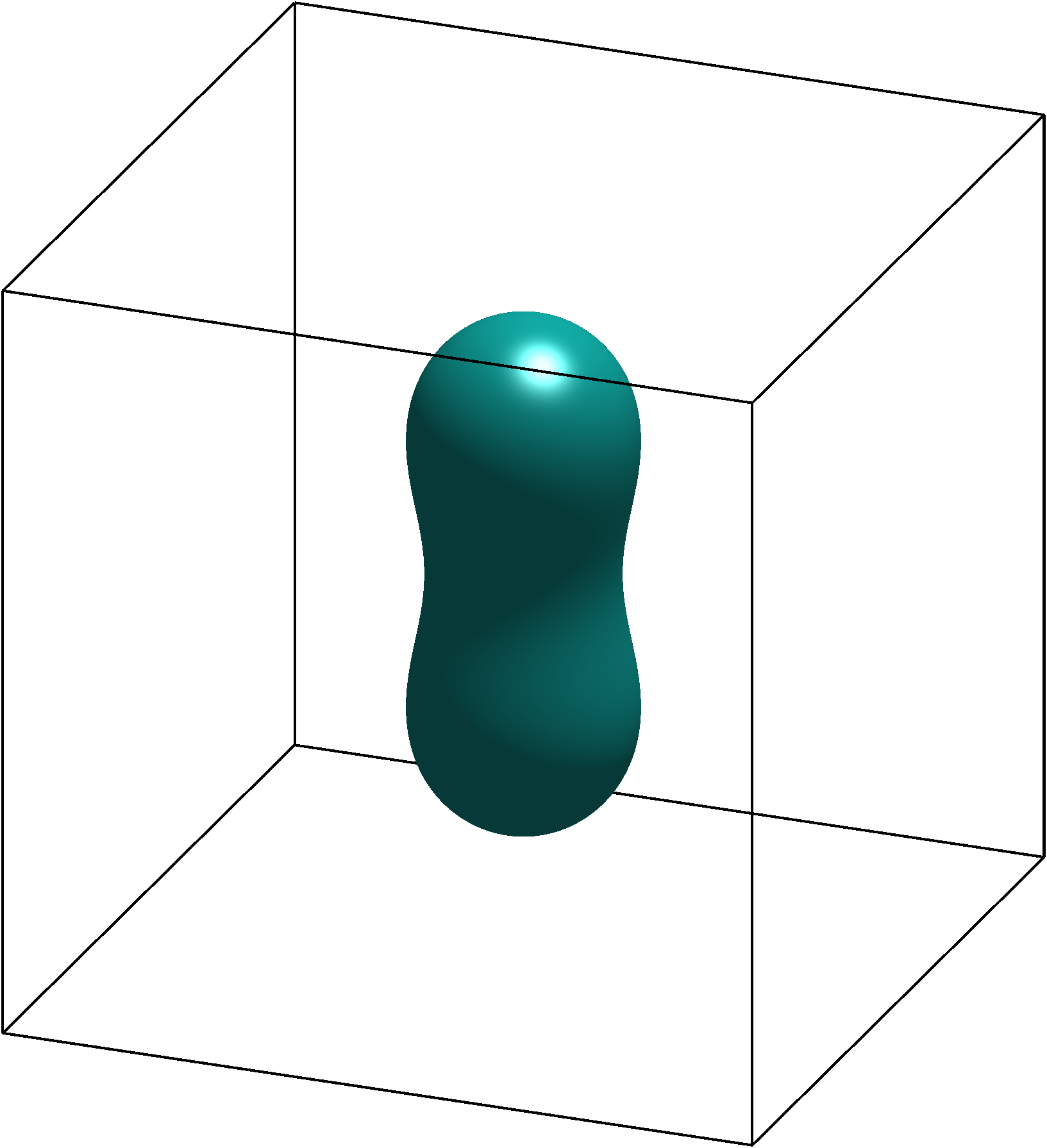}\\
        $t_{4}=0.8$  & $t_{5}=1.6$  & $t_{6}=3.2$ & $t_{7}=10$\\
        \includegraphics[width=0.21\linewidth]{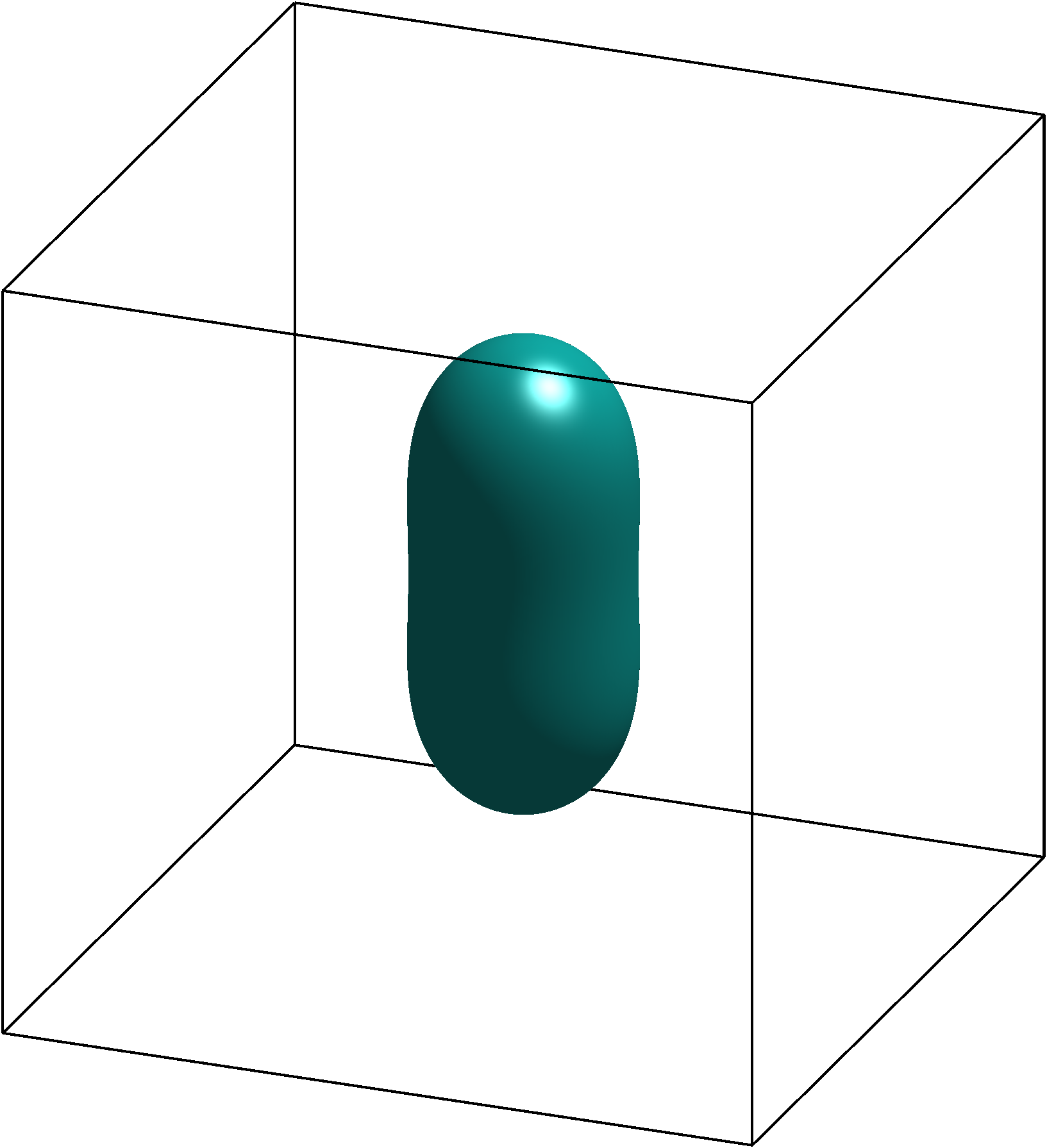}&
        \includegraphics[width=0.21\linewidth]{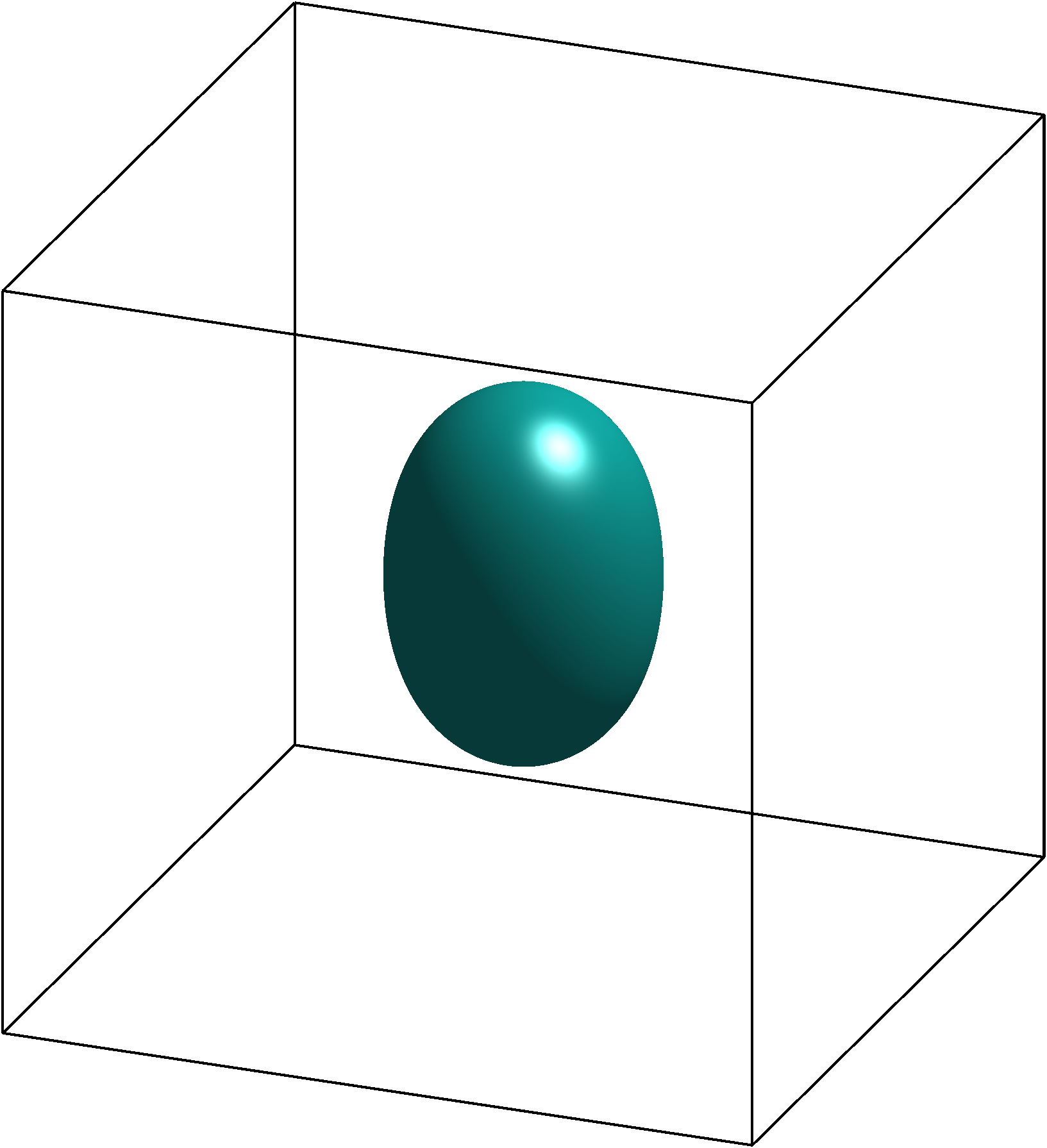}&
        \includegraphics[width=0.21\linewidth]{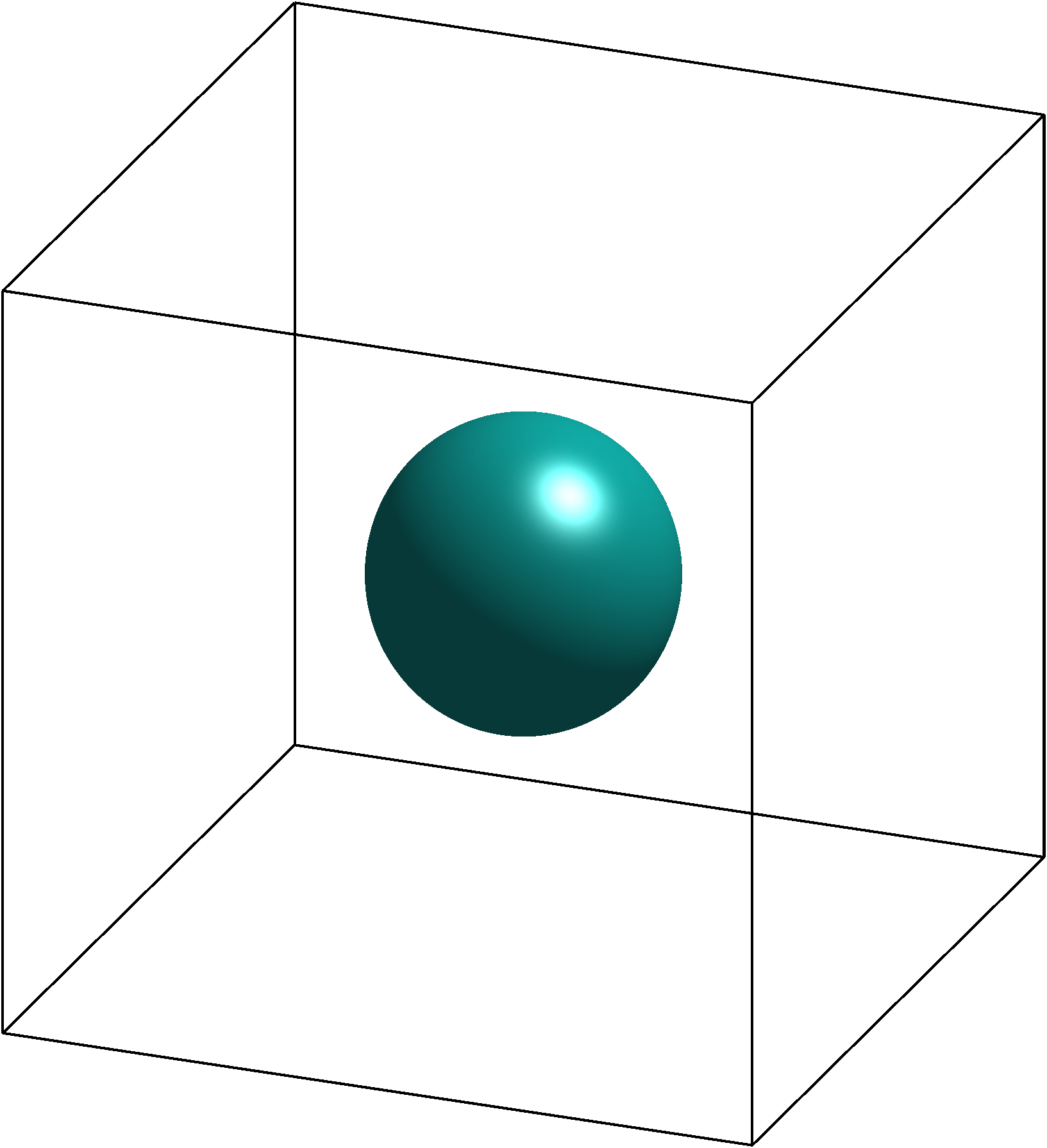}&
        \includegraphics[width=0.21\linewidth]{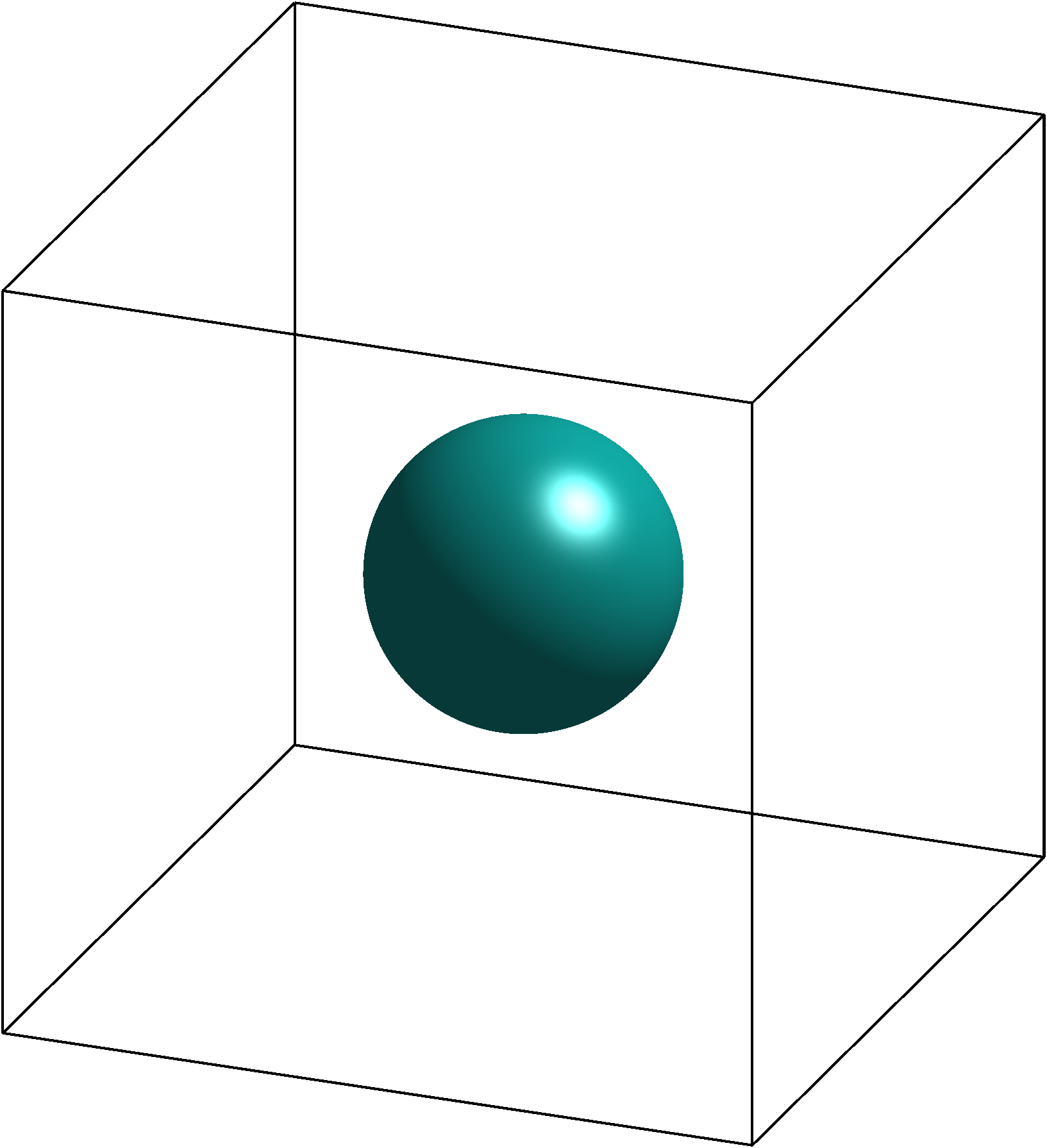}\\
    	\end{tabular}
    \end{center}
    \caption{\it Snapshots of the zero-isocontour of the phase function  $\phi$ show the coalescence of two drops at different time instants as indicated. See Table \ref{tab5: Cahn-Hilliard_GPU} for the computational time.}
    \label{fig3:CH_snap}
\end{figure}

Owing to the mass conservation and energy dissipation of the system \eqref{CH-sys}, the energy $E(\phi)$ first decreases before stabilizing at a constant value, as shown in Figure \ref{fig2:CH_energy}. The coalescence dynamics of the two droplets is illustrated through a series of temporal snapshots in Figure \ref{fig3:CH_snap}. These snapshots capture the evolving interfaces between the materials, visualized by the level set of $\phi=0$.
In Table \ref{tab5: Cahn-Hilliard_GPU}, we enumerate the online computational costs associated with various total DoFs. As explained above, the online computational time at each time step is less than solving two Poisson equations.
\begin{table}[ht!]
    \centering
    \begin{tabular}{|c|c|c|c|c|c|c|}
    \hline
         Total DoFs & $501^3$ & $551^3$ & $601^3$ & $651^3$ & $701^3$ & $751^3$\\
    \hline
        Total time & $8.46$E$2$ & $1.29$E$3$ & $1.79$E$3$ & $2.53$E$3$ & $3.10$E$3$ & $4.02$E$4$\\
    \hline
         Time for each time step & $8.46$E-$2$ & $1.29$E-$1$ & $1.79$E-$1$ & $2.53$E-$1$ & $3.10$E-$1$ & $4.02$E-$1$\\
    \hline
    \hline
         Total DoFs & $801^3$ & $851^3$ & $901^3$ & $951^3$ & $1001^3$ & $1051^3$\\
    \hline
         Total time & $5.41$E$3$ & $6.79$E$3$ & $8.84$E$3$ & $1.03$E$4$ & $1.27$E$4$ & -\\
    \hline
         Time for each time step & $5.41$E-$1$ & $6.79$E-$1$ & $8.84$E-$1$ & $1.03$E$0$ & $1.27$E$0$ & -\\
    \hline
    \end{tabular}
    \caption{\it Online computational time in seconds for $Q^5$ SEM in the BDF2 scheme \eqref{chscheme-2} solving the Cahn--Hilliard equation with 10,000 time steps for computing the solution at $T=10$ on Nivida A100. The {\it Total time} represents the online computational time for 10,000 time steps, and the {\it Time for each time step} is average online computational time per time step.}
    \label{tab5: Cahn-Hilliard_GPU}
\end{table}

\subsection{Comparison with implementation in Python}

For implementing \eqref{3D-Poisson}  on both CPU and GPU, similar to the implementation in MATLAB shown in Table \ref{tab:3dcode}, \eqref{3D-Poisson} can be  efficiently implemented  using the function {\it jax.numpy.einsum} in the Python package JAX  as shown in Table \ref{code:SEM-python}. 
\begin{table}[htbp]
    \centering
    \begin{pythonLines}
    u = jnp.einsum('ijk,kl->ijl',f,invTz.transpose())
    u = jnp.einsum('ijk,jl->ilk',u,invTy.transpose())
    u = jnp.einsum('li,ijk->ljk',invTx,u)
    u = u/Eig3D
    u = jnp.einsum('ijk,kl->ijl',u,Tz.transpose())
    u = jnp.einsum('ijk,jl->ilk',u,Ty.transpose())
    u = jnp.einsum('li,ijk->ljk',Tx,u)
    \end{pythonLines}
    \caption{\it The Python script of implementing (14) on both CPU and
GPU where \texttt{jnp} means \texttt{jax.numpy}.}
    \label{code:SEM-python}
\end{table}

Since both Python and MATLAB allow similar simple implementations of \eqref{3D-Poisson}  on GPU, it is interesting to compare them. 
We compare the performance of MATLAB with Python under double precision, as well as single precision, which often depends on specific hardware and their driver versions. 

 In Table \ref{tab8:poisson_python}, we list the online computational time comparison of similar implementations in MATLAB and Python on A100 for solving a 3D Poisson equation $200$ times. As we can see in Table \ref{tab8:poisson_python}, for double precision computation and problems with size smaller than $1000^3$, there is no significant difference in the online computational time between MATLAB and Python on GPUs. However, on A100 with 80G memory, MATLAB allows a problem size as large as $1250^3$, for which Python can handle only with single precision computation. \textcolor{blue}{It is noteworthy that the performance of single precision computation in Python can be significantly affected by the use of the TF32 (TensorFloat-32) format on NVIDIA GPUs. TF32 is a lower-precision format that requires less memory bandwidth and compute resources compared to the traditional FP32 (Float-32) format, leading to faster computation at the cost of some loss in precision.}

\begin{table}[htbp]
    \centering
    \begin{tabular}{|c|c|c|c|c|c|}
    \hline
     \multirow{2}{*}{Total DoFs}& \multicolumn{3}{c|}{Python(JAX)} &\multicolumn{2}{c|}{MATLAB}\\
     \cline{2-6} 
          & Single(TF32) & Single(FP32) & Double  & Single  & Double\\
    \hline
          $200^3$ & $4.80$E-$1$ & $6.88$E-$1$ & $6.70$E-$1$ & $4.72$E-$1$ & $5.20$E-$1$ \\
    \hline
          $250^3$ & $5.85$E-$1$ & $9.92$E-$1$ & $1.17$E$0$ & $7.32$E-$1$ & $9.11$E-$1$\\
    \hline
          $300^3$ & $7.43$E-$1$ & $1.67$E$0$ & $2.10$E$0$ & $1.54$E$0$& $1.80$E$0$ \\
    \hline
          $350^3$ & $1.39$E$0$ & $2.80$E$0$ & $3.47$E$0$ & $2.56$E$0$& $3.04$E$0$ \\
    \hline
          $400^3$ & $1.51$E$0$ & $4.56$E$0$ & $5.47$E$0$ & $4.63$E$0$& $5.07$E$0$ \\
    \hline
          $450^3$ & $2.89$E$0$ & $7.24$E$0$ & $8.74$E$0$ & $7.05$E$0$& $8.09$E$0$ \\
    \hline
          $500^3$ & $2.92$E$0$ & $9.48$E$0$ & $1.20$E$1$ & $8.92$E$0$& $1.07$E$1$ \\
    \hline
          $550^3$ & $6.07$E$0$ & $1.41$E$1$ & $1.82$E$1$ & $1.45$E$1$& $1.65$E$1$ \\
    \hline
          $600^3$ & $5.41$E$0$ & $1.99$E$1$ & $2.40$E$1$ & $1.97$E$1$& $2.28$E$1$ \\
    \hline
          $650^3$ & $1.02$E$1$ & $2.80$E$1$ & $3.35$E$1$ & $3.02$E$1$& $3.32$E$1$ \\
    \hline
          $700^3$ & $9.76$E$0$ & $3.41$E$1$ & $4.24$E$1$ & $3.55$E$1$& $4.04$E$1$ \\
    \hline
          $750^3$ & $1.66$E$1$ & $4.48$E$1$ & $5.66$E$1$ & $4.57$E$1$& $5.22$E$1$ \\
    \hline
          $800^3$ & $1.45$E$1$ & $6.12$E$1$ & $7.14$E$1$ & $6.32$E$1$& $6.89$E$1$ \\
    \hline
          $850^3$ & $2.66$E$1$ & $7.92$E$1$ & $9.28$E$1$ & $8.01$E$1$& $9.05$E$1$ \\
    \hline
          $900^3$ & $2.37$E$1$ & $1.00$E$2$ & $1.15$E$2$ & $1.05$E$2$& $1.15$E$2$ \\
    \hline
          $950^3$ & $4.10$E$1$ & $1.17$E$2$ & $1.37$E$2$ & $1.20$E$2$& $1.35$E$2$ \\
    \hline
          $1000^3$ & $3.17$E$1$ & $1.40$E$2$ & - & $1.40$E$2$ & $1.60$E$2$ \\
    \hline
          $1050^3$ & $6.12$E$1$ & $1.81$E$2$ & - & $1.87$E$2$ & $2.04$E$2$ \\
    \hline
          $1100^3$ & $4.77$E$1$ & $2.13$E$2$ & - & $2.16$E$2$ & $2.51$E$2$ \\
    \hline
          $1150^3$ & $7.96$E$1$ & $2.37$E$2$ & - & $2.39$E$2$ & $2.91$E$2$ \\
    \hline
          $1200^3$ & $6.30$E$1$ & $2.94$E$2$ & - & $2.95$E$2$ & $3.34$E$2$ \\
    \hline
          $1250^3$ & $1.13$E$2$ & $3.36$E$2$ & - & $3.45$E$2$ & $4.13$E$2$ \\
    \hline
    \end{tabular}
    \caption{\it Online computational time of single precision and double precison on one Nvidia A100 80G GPU card, for $Q^5$ SEM for solving a 3D Poisson equation $200$ times. The time unit is second. For double precision computation in Python on A100, an out-of-memory error will emerge for problems with size larger than $950^3$.   }
    \label{tab8:poisson_python}
\end{table}

As shown in Table \ref{tab8:poisson_python}, for larger problems such as one billion DoFs, the fastest implementation on GPU is Python with single precision computation, which might be suitable for some practical simulations. \textcolor{blue}{When using the default TF32 setting on A100, the accuracy of the $Q^5$ spectral-element method deteriorates for larger problem sizes, with the order of accuracy dropping below the expected value. However, by setting the environment variable \texttt{NVIDIA\_TF32\_OVERRIDE=0} to disable TF32 precision for certain operations, such as matrix multiplication and \texttt{jax.numpy.einsum}, the accuracy can be improved, and the expected order of accuracy is maintained.}{ To investigate the impact of TF32 on the accuracy of single precision computation in Python, we conducted additional tests with and without the TF32 format. Table \ref{tab10:accuracy test python} presents the accuracy results for the $Q^5$ spectral-element method under single precision on A100 with both TF32 and FP32 format.} In general, the implementation for high order SEM with \textcolor{blue}{TF32} single precision is not  robust  on A100, e.g., computation with SEM for the problem in Figure \ref{fig3:CH_snap} might blow up.

\begin{table}[htbp]
    \centering
    \begin{tabular}{|c|c c c|c c c|}
    \hline
    \multicolumn{7}{|c|}{$Q^5$ spectral-element method (\textcolor{blue}{TF32} single precision)}\\
    \hline
    \multirow{2}{*}{FEM Mesh} &  \multicolumn{3}{c|}{Dirichlet boundary} &  \multicolumn{3}{c|}{Neumann boundary}\\
     & {Total DoFs}  & $\ell^2$ error & order & {Total DoFs}  & $\ell^2$ error & order\\
    \hline
    $2^3$ &$9^3$ &  $2.27$E-$1$ & - & $11^3$ & $4.82$E-$1$ & - \\
    \hline
    $4^3$ & $19^3$ &  $3.92$E-$3$ & $5.86$ & $21^3$ & $6.60$E-$3$ & $6.19$\\
    \hline
    $8^3$ & $39^3$ &  $4.12$E-$5$ & $6.57$ & $41^3$ & $4.32$E-$5$ & $7.26$\\
    \hline
    $16^3$ & $79^3$ & $1.44$E-$3$ & -$5.13$ & $81^3$ & $2.46$E-$3$ & -$5.83$\\
    \hline
    $32^3$ & $159^3$ & $1.95$E-$3$ & -$0.44$ & $161^3$ & $2.73$E-$3$ & -$0.15$\\
    \hline
    \hline
    \multicolumn{7}{|c|}{$Q^5$ spectral-element method (\textcolor{blue}{FP32} single precision)}\\
    \hline
    \multirow{2}{*}{FEM Mesh} &  \multicolumn{3}{c|}{Dirichlet boundary} &  \multicolumn{3}{c|}{Neumann boundary}\\
     & {Total DoFs}  & $\ell^2$ error & order & {Total DoFs}  & $\ell^2$ error & order\\
    \hline
    $2^3$ &$9^3$ &  $2.27$E-$1$ & - & $11^3$ & $4.76$E-$1$ & - \\
    \hline
    $4^3$ & $19^3$ &  $3.91$E-$3$ & $5.86$ & $21^3$ & $5.49$E-$3$ & $6.44$\\
    \hline
    $8^3$ & $39^3$ &  $4.11$E-$5$ & $6.57$ & $41^3$ & $4.32$E-$5$ & $6.99$\\
    \hline
    $16^3$ & $79^3$ & $1.67$E-$6$ & $4.62$ & $81^3$ & $1.63$E-$6$ & $4.72$\\
    \hline
    $32^3$ & $159^3$ & $1.34$E-$6$ & $0.32$ & $161^3$ & $1.95$E-$6$ & -$0.26$\\
    \hline
    \end{tabular}
    \caption{\it Accuracy tests under \textcolor{blue}{TF32 and FP32} single precision in Python on Nvidia GPU A100 for the 3D Poisson equation \eqref{pde}
with $\alpha = 1$.  The actual accuracy of single precision computation depends very much on the hardware and version of  hardware 
 drivers. For Python, we implement the code under the environment JAX version 0.4.19 for Nvidia GPU A100, with Driver Version 535.86.10 and CUDA Version 12.2. See Table \ref{tab1:accuracy test dirichlet/neumann} for the results of MATLAB with double precision on Nvidia GPU A100.
    }
    \label{tab10:accuracy test python}
\end{table}

\begin{table}[htbp]
    \centering
    \begin{tabular}{|c|c c|c c|}
    \hline
    \multicolumn{5}{|c|}{FFT implementation on A100 for periodic boundary}\\
    \hline
    \multirow{2}{*}{Total DoFs} &  \multicolumn{2}{c|}{\textcolor{blue}{TF32} Single precision} &  \multicolumn{2}{c|}{Double precision}\\
     &  $\ell^2$ error & order & $\ell^2$ error & order\\
    \hline
    $10^3$ &  $5.00$E-$1$ & -  & $5.00$E-$1$ & - \\
    \hline
    $20^3$ &  $1.05$E-$1$ & $2.25$ & $1.05$E-$1$ & $2.25$\\
    \hline
    $40^3$ &  $2.53$E-$2$ & $2.06$ & $2.53$E-$2$ & $2.06$\\
    \hline
    $80^3$ &  $6.26$E-$3$ & $2.01$ & $6.26$E-$3$ & $2.01$\\
    \hline
    $160^3$ & $1.56$E-$3$ & $2.01$ & $1.56$E-$3$ & $2.00$\\
    \hline
    $320^3$ & $3.88$E-$4$ & $2.00$ & $3.90$E-$4$ & $2.00$\\
    \hline
    $640^3$ & $8.57$E-$5$ & $2.18$ & $9.75$E-$5$ & $2.00$\\
    \hline
    $900^3$ & $5.64$E-$5$ & $1.23$ & $4.93$E-$5$ & $2.00$\\
    \hline
    $1200^3$ & $9.32$E-$5$ & -$1.75$ & - & -\\
    \hline
    \end{tabular}
    \caption{\it Accuracy tests for second order finite difference with periodic boundary (FFT implementation) in Python on Nvidia GPU A100 for the 3D Poisson equation \eqref{pde} with $u^{*} = \sin(2\pi x)\sin(3\pi y)\sin(4\pi z)$ and $\alpha = 1$.}
    \label{tab11:accuracy test python}
\end{table}

\begin{figure}[ht!]
    \begin{center}
	   \begin{tabular}{cccc}
        $t_{0}=0$  & $t_{1}=0.1$  & $t_{2}=0.2$ & $t_{3}=0.4$\\
        \includegraphics[width=0.21\linewidth]{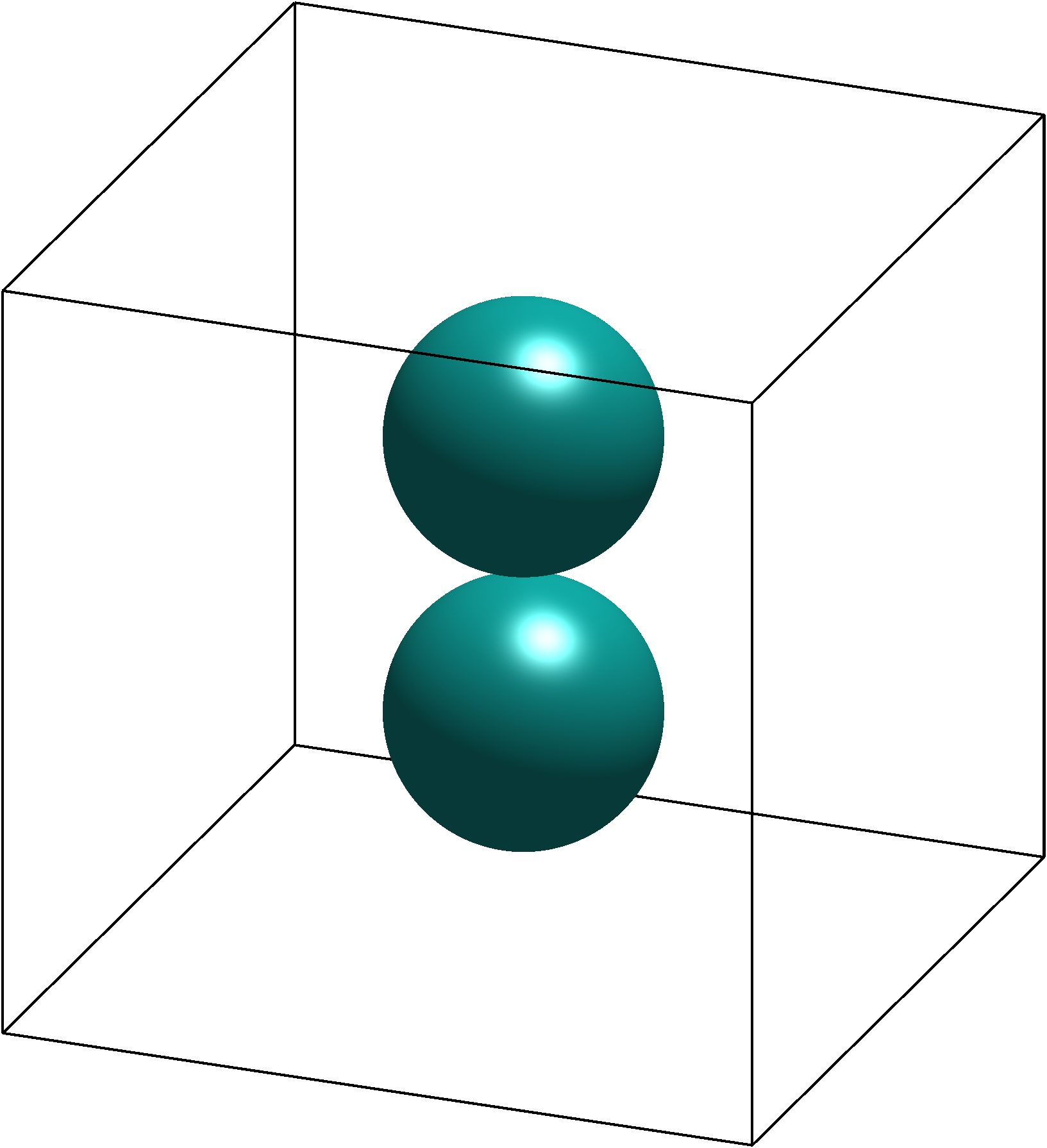}&
        \includegraphics[width=0.21\linewidth]{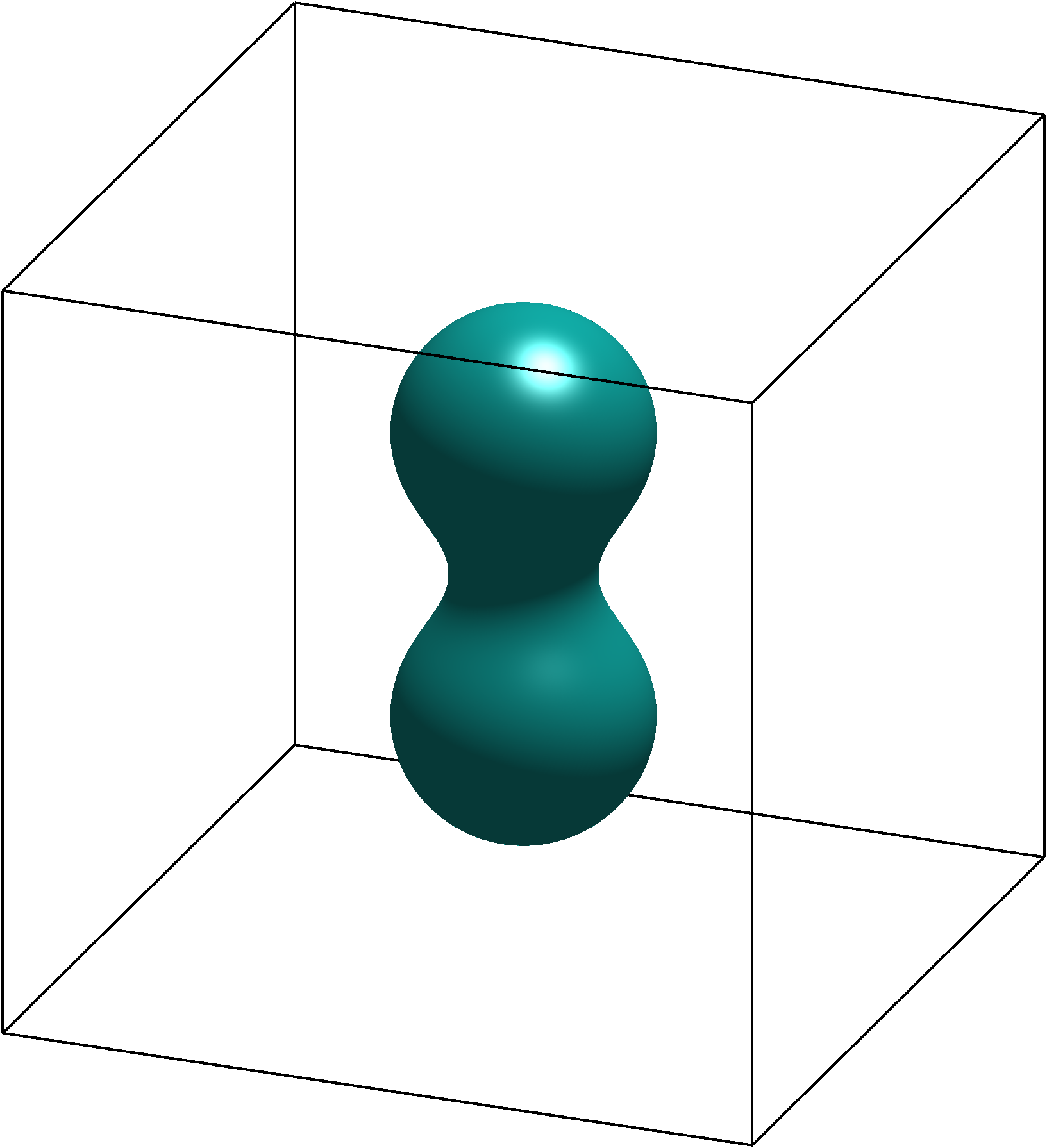}&
        \includegraphics[width=0.21\linewidth]{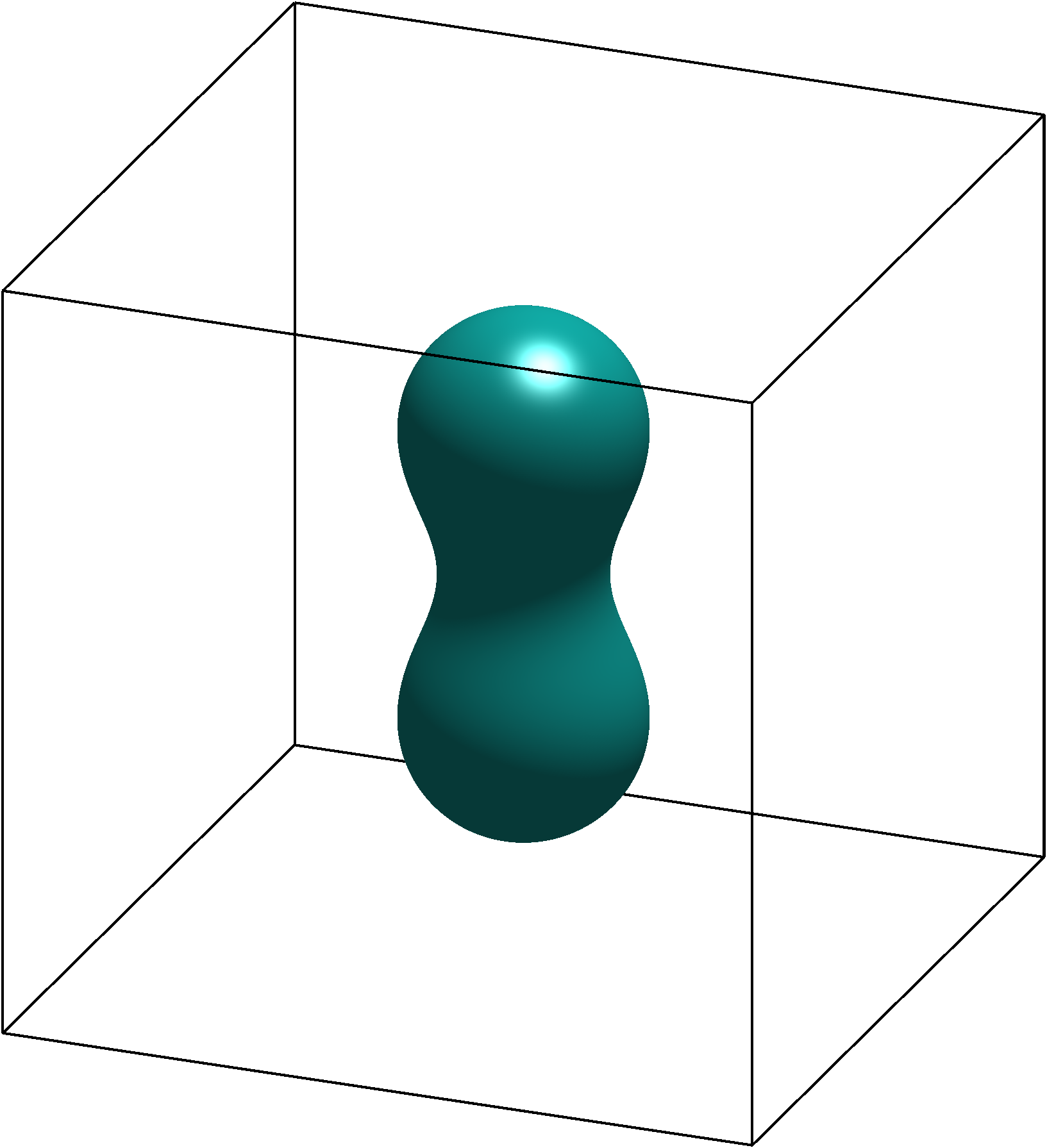}&
        \includegraphics[width=0.21\linewidth]{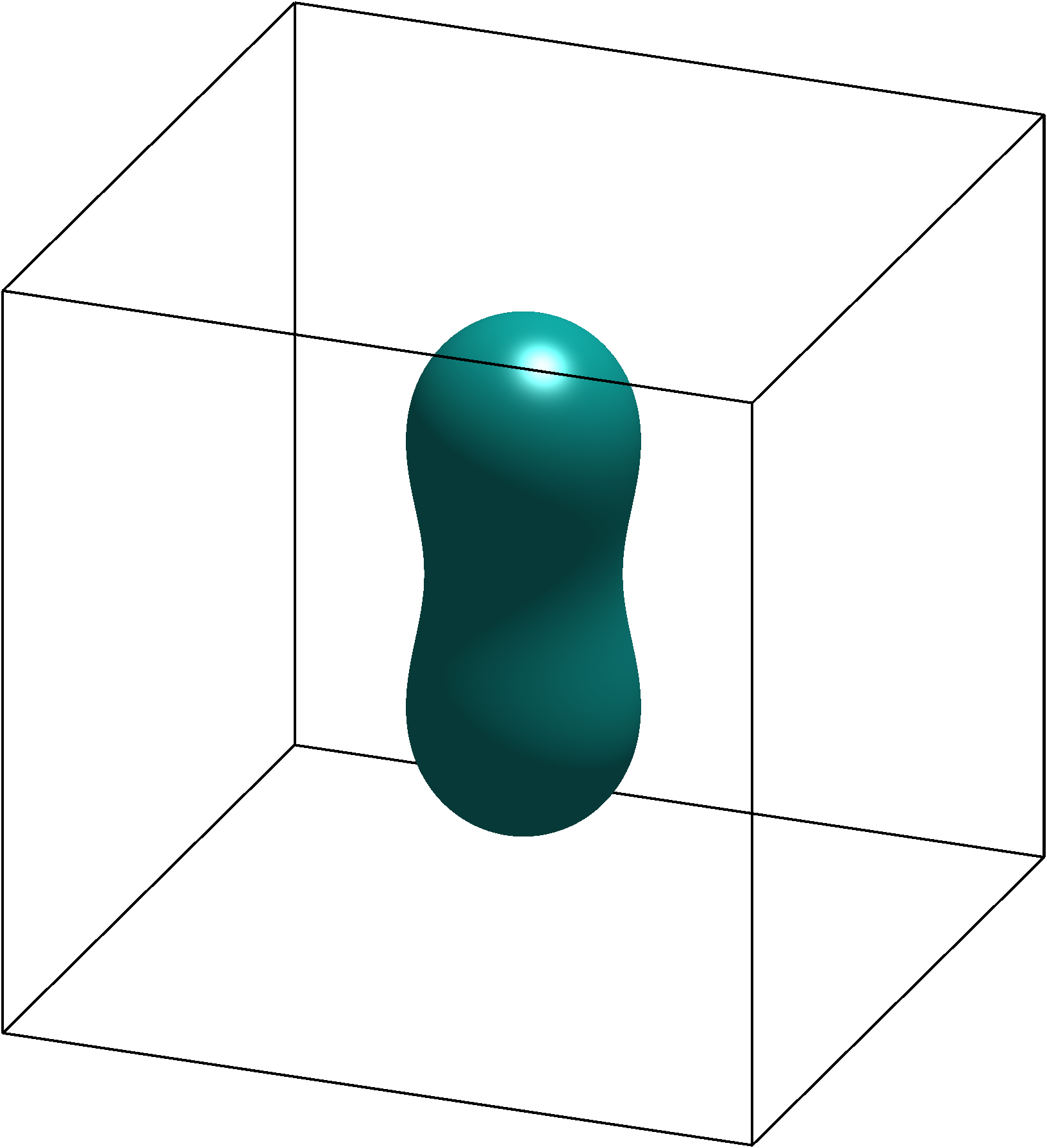}\\
        $t_{4}=0.8$  & $t_{5}=1.6$  & $t_{6}=3.2$ & $t_{7}=10$\\
        \includegraphics[width=0.21\linewidth]{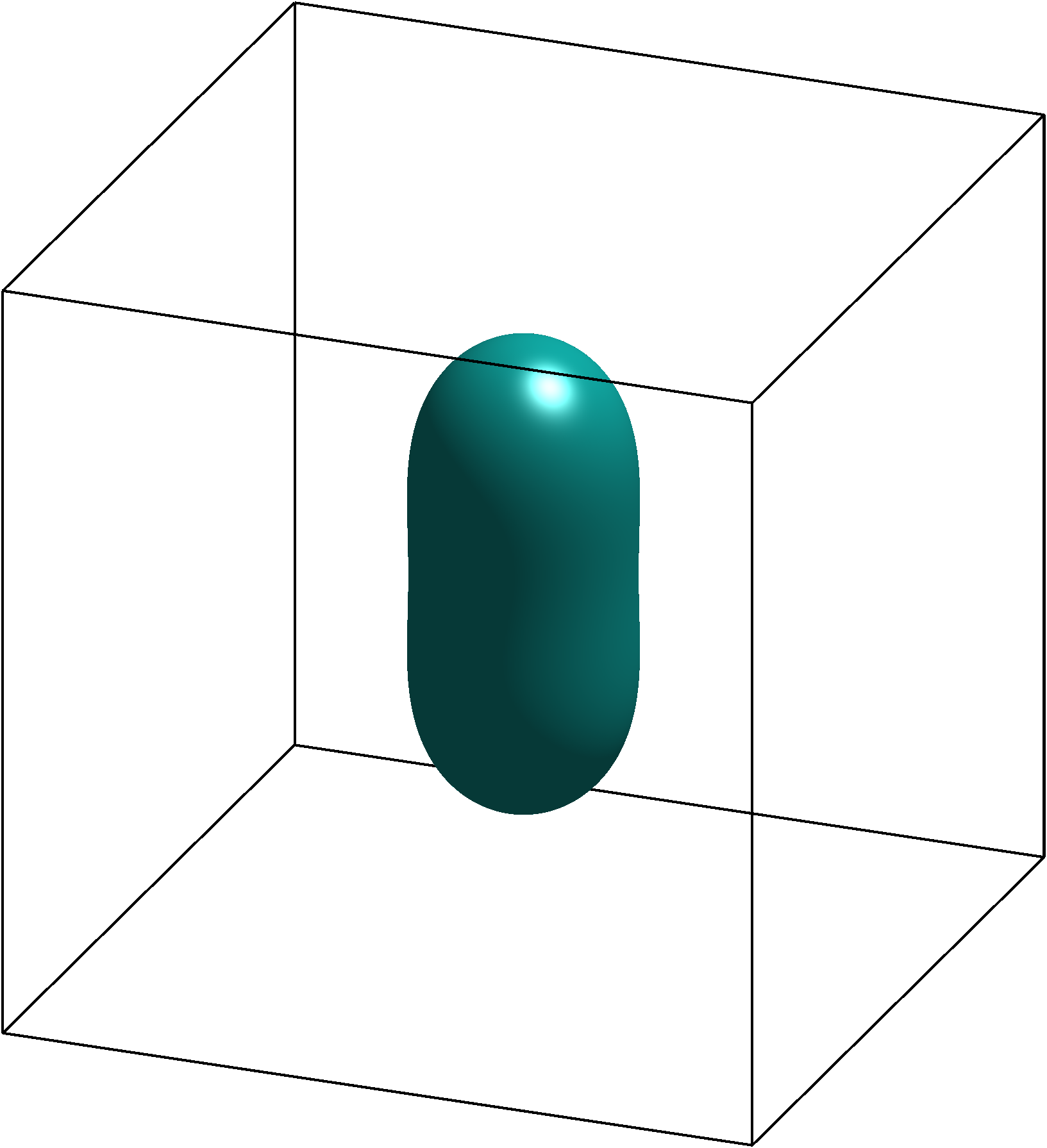}&
        \includegraphics[width=0.21\linewidth]{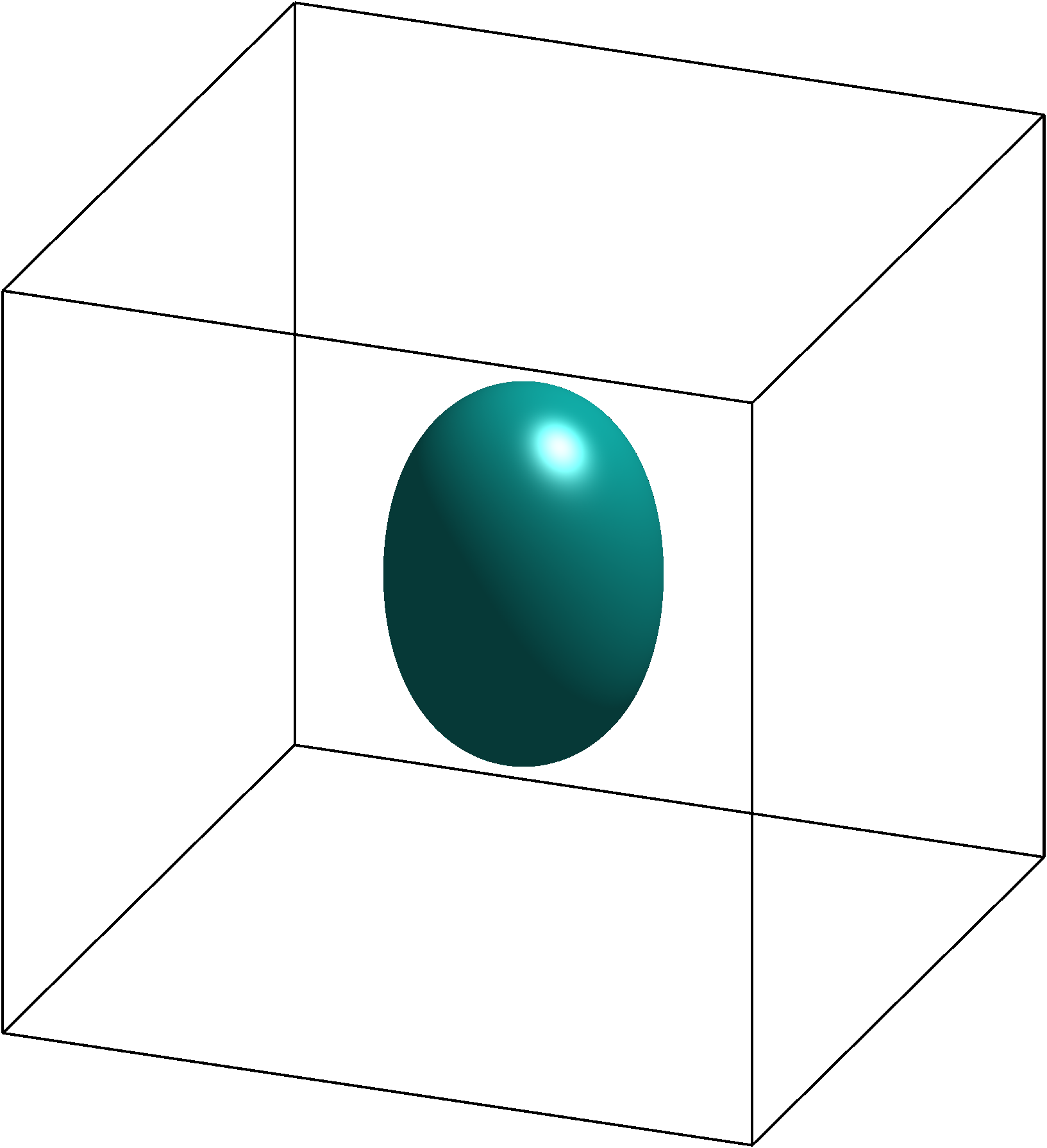}&
        \includegraphics[width=0.21\linewidth]{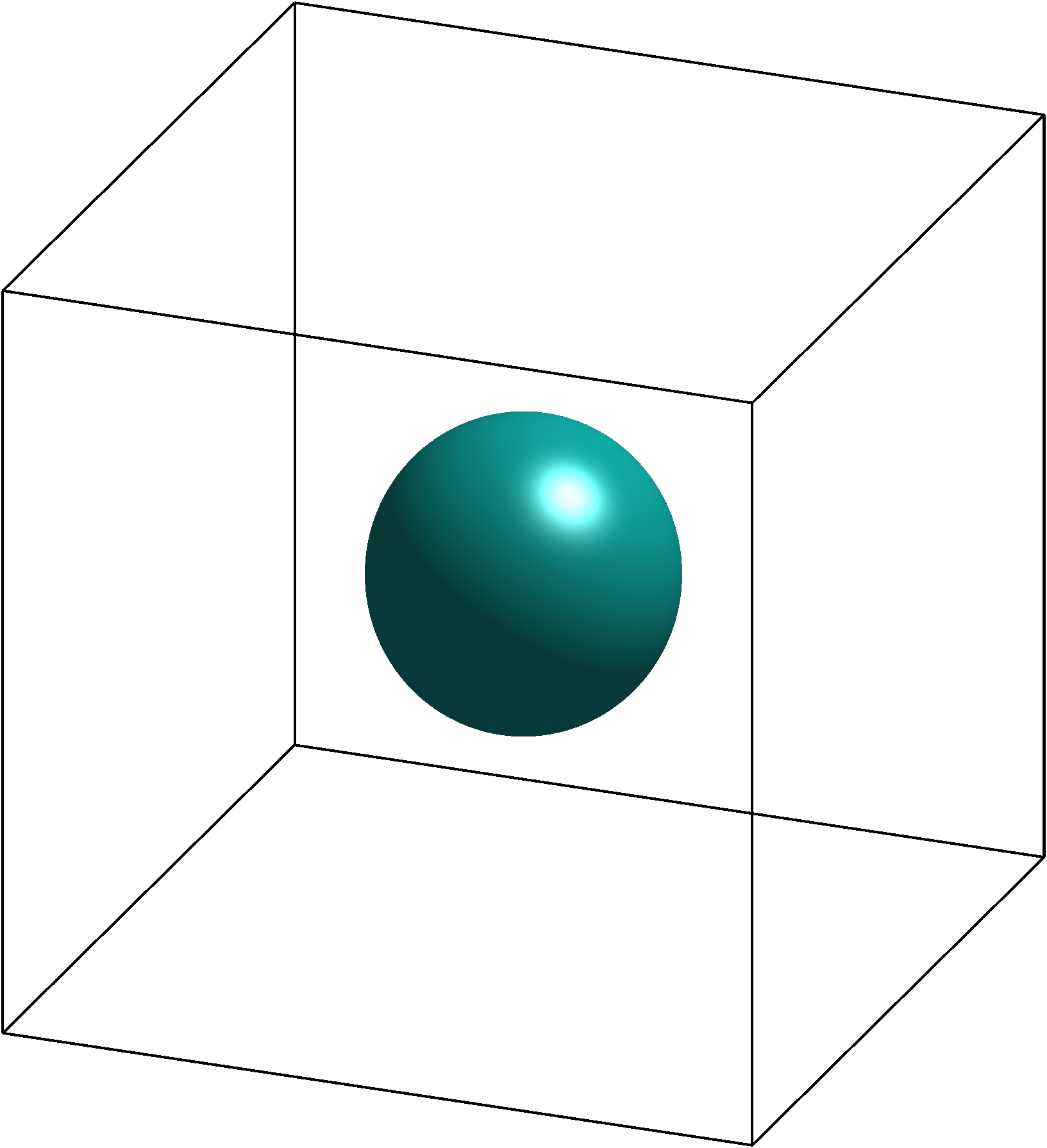}&
        \includegraphics[width=0.21\linewidth]{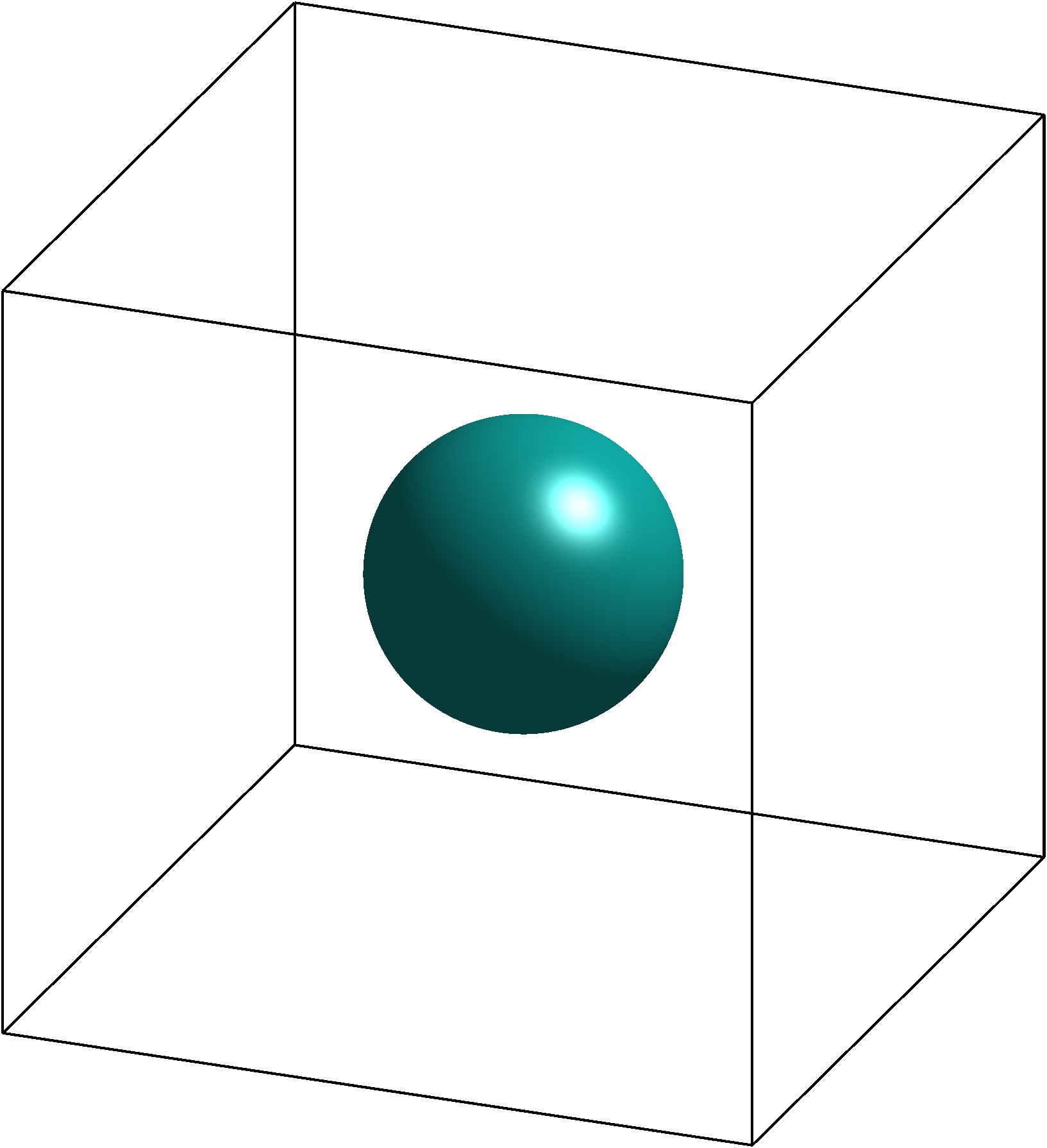}\\
    	\end{tabular}
    \end{center}
    \caption{\it Snapshots of the same problem in Figure \ref{fig3:CH_snap}, implemented by Python in  \textcolor{blue}{TF32} single precision on A100 for second order finite difference (i.e., $Q^1$ spectral-element method) by FFT with total DoFs $800^3$.}
    \label{fig5:CH_snap_python_FFT}
\end{figure}

On the other hand, the second order finite difference implemented in Python Jax with \textcolor{blue}{TF32} single precision computation \textcolor{blue}{is robust as suggested by Table \ref{tab11:accuracy test python}}. For periodic boundary conditions, the eigenvectors of second order finite difference (i.e., $Q^1$ spectral-element method) can be implemented by FFT as shown in Table \ref{code:FFT-python}. 
As a demonstration, 
we include the computation results for the Cahn-Hilliard equation of Python in  \textcolor{blue}{TF32} single precision on A100 in Figure \ref{fig5:CH_snap_python_FFT}, which is comparable to the double precision results on A100 in Figure \ref{fig3:CH_snap}.

\begin{table}[htbp]
    \centering
    \begin{pythonLines}
    u = jnp.fft.fftn(f)/Eig3D
    if alpha == 0:
        u[0,0,0] = 0.
    u = jnp.real(jnp.fft.ifftn(u))
    \end{pythonLines}
    \caption{\it The Python script for FFT implementation of a second order (i.e., $Q^1$ spectral-element method) for the Poisson equation with periodic boundary conditions on both CPU and GPU where \texttt{jnp} means \texttt{jax.numpy}.}
    \label{code:FFT-python}
\end{table}


\section{Concluding remarks}
\label{sec:remark}
In this paper, we have discussed a simple MATLAB 2023 implementation for accelerating high order methods on GPUs. For large enough 3D problems, a speed-up of at least 60 can be achieved on Nvidia A100. In particular, solving a 3D Poisson type equation with one billion DoFs costs only 0.8 second for $Q^k$ spectral-element method. As examples of applications, we  applied this fast solver to solve a linear (time-independent) Schr\"odinger equation and a nonlinear (time-dependent) Cahn-Hilliard equation in three-dimension. { 
We expect  the proposed simple implementation to have the same performance for  any problem with similar
tensor product structure, e.g., exponential time differencing and spectral fractional Laplacians. }
 
\section*{Data availability statements}
The authors declare that the data supporting the findings of this study are available within the paper and its supplementary information files.

\section*{Declarations}
 J. Shen's research was supported in part by NSFC  12371409, and X. Zhang's research was supported by NSF DMS-220815.
 The authors declare they have no financial interests.
\appendix 
\section*{Appendix}
\section{MATLAB scripts for a 3D Poisson equation}
\label{sec:appendix}
We provide a demonstration in MATLAB 2023  for $Q^k$ spectral-element method solving a Poisson equation in three dimensions, which involves three MATLAB scripts:
\begin{enumerate}
  \item   {\it Poisson3D.m} for solving the Poisson equation  on either CPU or GPU; 
  \item {\it SEGenerator1D.m} for generating stiffness and mass matrices in spectral element method;
\item  {\it LegendreD.m} for  Legendre and Jocaboi polynomials from \cite{hesthaven2007nodal}.
\end{enumerate}
Readers can easily reproduce the results in Section \ref{sec:poisson-accuracy} and Section \ref{sec:poisson-speedup} using these three MATLAB  scripts.

\newpage
\noindent{\it Poisson3D.m}:
\vspace{-0.5cm}
\begin{lstlisting}[numbers=none]
% Solving the Poisson equation by Q5 SEM with Neumann b.c.
if gpuDeviceCount('available')<1; Param.device='cpu';
else;Param.device='gpu'; Param.deviceID=1;end % ID=1,2,3,...
Np=5; Param.Np=Np; % polynomial degree Q5
Ncellx=40; Ncelly=40; Ncellz=40; % finite element cell number
% total number of unknowns in each direction
nx=Ncellx*Np+1; ny=Ncelly*Np+1; nz=Ncellz*Np+1;
% the domain is [-Lx, Lx]*[-Ly, Ly]*[-Lz, Lz]
Lx=1; Ly=1; Lz=1; cx=pi; cy=2*pi; cz=3*pi; alpha=1;
Param.Ncellx=Ncellx; Param.Ncelly=Ncelly; Param.Ncellz=Ncellz;
Param.nx = nx; Param.ny = ny; Param.nz = nz;
fprintf('3D Poisson with total DoFs %d by %d by %d \n',nx,ny,nz);
fprintf('Laplacian is Q%d spectral element method \n', Np);  
[x,ex,Tx,eigx]=SEGenerator1D('x',Lx,Param);
[y,ey,Ty,eigy]=SEGenerator1D('y',Ly,Param);
[z,ez,Tz,eigz]=SEGenerator1D('z',Lz,Param);
% a smooth solution
u1x=cos(cx*x); u2x=power(1-power(x,2),3); 
du2x=30*power(x,4)-36*power(x,2)+6;
u1y=cos(cy*y); u2y=power(1-power(y,2),2); 
du2y=4-12*power(y,2);
u1z=cos(cz*z); u2z=power(1-power(z,2),4); 
du2z=(8-56*power(z,2)).*power(1-power(z,2),2);
uexact=squeeze(tensorprod(u1x*u1y',u1z)+tensorprod(u2x*u2y',u2z));
f=(cx*cx+cy*cy+cz*cz)*squeeze(tensorprod(u1x*u1y',u1z))+...
  squeeze(tensorprod(du2x*u2y',u2z)+tensorprod(u2x*du2y',u2z)+...
  tensorprod(u2x*u2y',du2z))+alpha*uexact;
TxInv=pinv(Tx); TyInv=pinv(Ty); TzInv=pinv(Tz);
if strcmp(Param.device,'gpu');Device=gpuDevice(Param.deviceID);
   fprintf('GPU computation: starting to load matrices/data \n');
   Tx=gpuArray(Tx); Ty=gpuArray(Ty); Tz=gpuArray(Tz);
   eigx=gpuArray(eigx);eigy=gpuArray(eigy);eigz=gpuArray(eigz);
   ex=gpuArray(ex); ey=gpuArray(ey); ez=gpuArray(ez); f=gpuArray(f); 
   TxInv=gpuArray(TxInv);TyInv=gpuArray(TyInv);TzInv=gpuArray(TzInv);
end 
Lambda3D=squeeze(tensorprod(eigx,ey*ez')+tensorprod(ex,eigy*ez')...
        +tensorprod(ex,ey*eigz'));
if strcmp(Param.device,'gpu'); wait(Device);
    fprintf('GPU loading finished and computing started \n');
end
tic; % online computation
u = tensorprod(f,TzInv',3,1);
u = pagemtimes(u,TyInv');
u = squeeze(tensorprod(TxInv,u,2,1));
u = u./(Lambda3D + alpha);
u = tensorprod(u,Tz',3,1);
u = pagemtimes(u,Ty');
u = squeeze(tensorprod(Tx,u,2,1));
if strcmp(Param.device,'gpu');wait(Device);end;time=toc;err=u-uexact;
fprintf('The ell infinity norm error is %d \n',norm(err(:),inf));
if strcmp(Param.device,'gpu')
   fprintf('The online GPU computation time is %d \n', time);
else
   fprintf('The online CPU computation time is %d \n', time);
end
\end{lstlisting}

 \newpage
 \noindent {\it SEGenerator1D.m}:
\begin{lstlisting}[numbers=none]
function [varargout] = SEGenerator1D(direction,L,Param)
% generate 1D spectral element with Neumann B.C.
switch direction
    case 'x'
        N=Param.Np; Ncell=Param.Ncellx; n=Param.nx;
    case 'y'
        N=Param.Np; Ncell=Param.Ncelly; n=Param.ny;
    case 'z'
        N=Param.Np; Ncell=Param.Ncellz; n=Param.nz;
end
% generate the mesh with Ncell intervals with domain [Left, Right]
[D,r,w] = LegendreD(N); Left = -L; Right = L;
Length = Right - Left; dx = Length/Ncell;
for j = 1:Ncell
    cellLeft = Left+dx*(j-1);
    localPoints = cellLeft+dx/2+r*dx/2;
    if (j==1)
    x = localPoints;
    else
    x = [x;localPoints(2:end)];
    end
end        
SLocal = D'*diag(w)*D; % local stiffness matrix for each element       
S=[]; M=[];
for j = 1:Ncell % global stiffness and lumped mass matrices
    S = blkdiag(SLocal,S); M = blkdiag(diag(w),M); 
end
% Next step: "glue" the cells
Np = N+1; % number of points in each cell
Glue = sparse(zeros(Ncell*Np-Ncell+1, Ncell*Np));
for j = 1:Ncell
    rowStart=(j-1)*Np+2-j; rowEnd=rowStart+Np-1;
    colStart=(j-1)*Np+1; colEnd=colStart+Np-1;
    Glue(rowStart:rowEnd,colStart:colEnd)=speye(Np); 
end
S=Glue*S*Glue'; M=Glue*M*Glue'; H=diag(1./diag(M))*S; 
ex=ones(n,1); MHalfInv=diag(1./sqrt(diag(M)));
S1=MHalfInv*S*MHalfInv; S1=(S1+S1')/2;
[U,d]=eig(S1,'vector'); [lambda,indexSort]=sort(d);
T=U(:,indexSort); h=dx/2; lambda=lambda/(h*h); 
S1=sparse(S1/(h*h)); M=sparse(M);
% after this step, T is the eigenvector of H
T = MHalfInv*T; H=full(H/(h*h)); S=S/h; M=full(M*h);
varargout{1}=x; varargout{2}=ex; varargout{3}=T; varargout{4}=lambda;
end
\end{lstlisting}
\newpage

 \noindent{\it LegendreD.m}:
\begin{lstlisting}[numbers=none]
function [D,r,w] = LegendreD(N)
    Np = N+1; r = JacobiGL(0,0,N);
    w = (2*N+1)/(N*N+N)./power(JacobiP(r,0,0,N),2);   
    Distance = r*ones(1,N+1)-ones(N+1,1)*r'+eye(N+1);
    omega = prod(Distance,2);
    D = diag(omega)*(1./Distance)*diag(1./omega); 
    D(1:Np+1:end) = 0; D(1:Np+1:end) = -sum(D,2);
end
function [x] = JacobiGL(alpha,beta,N)
    x = zeros(N+1,1);
    if (N==1); x(1)=-1.0; x(2)=1.0; return; end
    [xint,temp] = JacobiGQ(alpha+1,beta+1,N-2);
    x = [-1, xint', 1]'; return;
end
function [x,w] = JacobiGQ(alpha,beta,N)
    if (N==0) 
        x(1) = -(alpha-beta)/(alpha+beta+2); w(1) = 2; return; 
    end
    h1 = 2*(0:N)+alpha+beta;
    J = diag(-1/2*(alpha*alpha-beta*beta)./(h1+2)./h1) + ...
        diag(2./(h1(1:N)+2).*sqrt((1:N).*((1:N)+alpha+beta).*...
        ((1:N)+alpha).*((1:N)+beta)./(h1(1:N)+1)./(h1(1:N)+3)),1);
    if (alpha+beta<10*eps); J(1,1)=0.0; end
    J = J + J'; [V,D] = eig(J); x = diag(D);
    w = power(V(1,:)',2)*power(2,alpha+beta+1)/(alpha+beta+1)*...
        gamma(alpha+1)*gamma(beta+1)/gamma(alpha+beta+1);
end
function [P] = JacobiP(x,alpha,beta,N)
    xp = x; dims = size(xp);
    if (dims(2)==1); xp = xp'; end    
    PL = zeros(N+1,length(xp));   
    gamma0 = power(2,alpha+beta+1)/(alpha+beta+1)*gamma(alpha+1)*...
        gamma(beta+1)/gamma(alpha+beta+1);
    PL(1,:) = 1.0/sqrt(gamma0);
    if (N==0); P = PL'; return; end
    gamma1 = (alpha+1)*(beta+1)/(alpha+beta+3)*gamma0;
    PL(2,:) = ((alpha+beta+2)*xp/2 + (alpha-beta)/2)/sqrt(gamma1);
    if (N==1); P = PL(N+1,:)'; return; end
    aold = 2/(2+alpha+beta)*sqrt((alpha+1)*(beta+1)/(alpha+beta+3));    
    for i = 1:N-1
      h1 = 2*i+alpha+beta;
      anew = 2/(h1+2)*sqrt( (i+1)*(i+1+alpha+beta)*(i+1+alpha)*...
          (i+1+beta)/(h1+1)/(h1+3));
      bnew = - (alpha*alpha-beta*beta)/h1/(h1+2);
      PL(i+2,:) = 1/anew*( -aold*PL(i,:) + (xp-bnew).*PL(i+1,:));
      aold = anew;
    end  
    P = PL(N+1,:)';
end
\end{lstlisting}

\bibliographystyle{plain} 
\bibliography{reference}
\end{document}